\documentclass[10pt,a4paper,leqno]{article}
\usepackage[T1]{fontenc}
\usepackage[latin1]{inputenc}
\usepackage{float}
\usepackage{graphics}
\usepackage{amssymb}
\usepackage{float}
\usepackage{makeidx}
\usepackage{graphicx}
\usepackage{amsmath} 
\usepackage{pstricks}

\catcode`\é=\active \def é{\' e}

\newfont{\Blackboard}{msbm10 scaled 1200}

\newfont{\roma}{cmr10 scaled 1200}
 
\def \Z {{\mathbb{Z}}}
\def \R {{\mathbb{R}}}
\def \N {{\mathbb{N}}}
\def \C {{\mathbb{C}}}
\def \Q {{\mathbb{Q}}}

\newcommand {\nc}   {\newcommand}

\def\dfrac{\displaystyle \frac }
\def\edc{\end{document}}

\nc {\be}   {\begin{equation}} \nc {\ee}   {\end{equation}} 
\nc {\beq}  {\begin{eqnarray}} \nc {\eeq}  {\end{eqnarray}}
\nc {\beqs} {\begin{eqnarray*}} \nc {\eeqs} {\end{eqnarray*}}
\newcommand{\bthe}{\begin{theorem}}
\newcommand{\ethe}{\end{theorem}}
\newcommand{\brk}{\begin{remark}}
\newcommand{\erk}{\end{remark}}
\newcommand{\bco}{\begin{corollary}}
\newcommand{\eco}{\end{corollary}}
\newcommand{\blem}{\begin{lemma}}
\newcommand{\elem}{\end{lemma}}
\newcommand{\bprop}{\begin{proposition}}
\newcommand{\eprop}{\end{proposition}}
\newcommand{\bex}{\begin{example}}
\newcommand{\eex}{\end{example}}
\newcommand{\fin}{\,\rule{1ex}{1ex}\,}
\newcommand {\proof} {\noindent \bf Proof. \rm}

\newcommand{\caD}{{\cal D}}

\newcommand{\caA}{{\cal A}}

\newcommand{\caL}{{\cal L}}

\newcommand{\caH}{{\cal H}}

\newcommand{\da}{\caD(\caA)}

\newtheorem{theorem}{Theorem}[section]
\newtheorem{lemma}[theorem]{Lemma}
\newtheorem{corollary}[theorem]{Corollary}
\newtheorem{remark}[theorem]{Remark}
\newtheorem{definition}[theorem]{Definition}
\newtheorem{proposition}[theorem]{Proposition}
\newtheorem{example}[theorem]{Example}
\title{Do Shape Memory Alloy cables restrain the vibrations of girder bridges? - a mathematical point of view}

\author{V. Régnier\thanks{Laboratoire de Mathématiques
et ses Applications de Valenciennes, FR CNRS 2956, Institut des Sciences et Techniques de
Valenciennes, Université Polytechnique des Hauts de France, Le Mont
Houy, 59313 VALENCIENNES Cedex 9, FRANCE, e-mail : virginie.regnier@uphf.fr}}

\begin{document}

\maketitle




\begin{abstract}
We study the energy decay of a damped Euler-Bernoulli beam which is subject to a pointwise feedback force representing a Shape Memory Alloy (SMA) cable. The problem we consider is that of \cite{LiuFu} but, for simplicity, our modelization does not take into account the additional stiffness term they considered. An explicit expression is given for the resolvent of the underlying operator as well as its eigenvalues and eigenfunctions. We show the exponential decay of the energy. The fastest decay rate is given by the supremum of the real part of the spectrum of the infinitesimal generator of the underlying semigroup since we prove the existence of a Riesz basis. To the question "Do Shape Memory Alloy cables restrain the vibrations of girder bridges?", the experiments in \cite{LiuFu} answer positively. Our study does not allow to give a definite answer yet. The only presence of these cables may not to be enough. Some physical parameters have to be chosen carefully. 
\end{abstract}



\maketitle

\section{Introduction}
\setcounter{equation}{0}

In this paper, we consider the same problem as in \cite{LiuFu} which is rewritten for simplicity:

\begin{eqnarray}
\label{eq1} 
\dfrac{\partial^2 u}{\partial t^2}(x,t) + a \dfrac{\partial^4 u}{\partial x^4}(x,t) + b \dfrac{\partial u}{\partial t}(x,t) + \alpha \dfrac{\partial u}{\partial t}(\xi,t) \delta_{\xi}  &=& 0  \hskip 1cm \mbox{for } (x,t) \;  \mbox{in } (0,1) \times (0,\infty) \\
u(0,t)=u(1,t)= \dfrac{\partial^2 u}{\partial x^2}(0,t) = \dfrac{\partial^2 u}{\partial x^2}(1,t) &=&0  \hskip 1cm \mbox{for } t \; \mbox{in } (0,\infty)\\
\label{BC}
u(x,0)=u_0(x), \; \; \dfrac{\partial u}{\partial t}(x,0)=u_1(x) & & \hskip 1.2cm \mbox{for } x \; \mbox{in } (0,1) 
\label{IC}
\end{eqnarray}

\noindent where $a$, $b$ and $\alpha$ are strictly positive constants. \\
\\
In this problem the function $u$ denotes the transverse displacement of the bridge deck represented by a beam with a uniform section, $\delta_{\xi}$ is the Dirac mass concentrated in the point $\xi \in (0;1)$.\\
\\
\noindent Denote by $\rho$ the mass density of the beam, $A$ the area of the cross section of the beam, $EI$ the second moment of area of the cross-section and $c$ the damping of the beam. Then (\ref{eq1})-(\ref{BC}) coincides with the system in \cite{LiuFu} with $a =\dfrac{(EI)}{\rho A}$, $b= \dfrac{c}{\rho A}$, $\alpha= \dfrac{c_{SMA}}{\rho A}$ and $k_{SMA} = 0$.  \\
The values $k_{SMA}$ and $c_{SMA}$ are respectively the equivalent stiffness and equivalent damping of a damper which represents the Shape Memory Alloy (SMA) cable at the position $x$. They are calculated by $(26)$ and $(27)$ of \cite{LiuFu}. \\
\\
In fact, the problem considered in \cite{LiuFu} is more precisely the following one with $\beta > 0$:

\begin{eqnarray}
\label{eq1bis} 
\dfrac{\partial^2 u}{\partial t^2}(x,t) + a \dfrac{\partial^4 u}{\partial x^4}(x,t) + b \dfrac{\partial u}{\partial t}(x,t) + \left( \alpha \dfrac{\partial u}{\partial t}(\xi,t) + \beta u(\xi,t) \right) \delta_{\xi}  &=& 0  \hskip 1cm \mbox{for } (x,t) \;  \mbox{in } (0,1) \times (0,\infty) \\
u(0,t)=u(1,t)= \dfrac{\partial^2 u}{\partial x^2}(0,t) = \dfrac{\partial^2 u}{\partial x^2}(1,t) &=&0  \hskip 1cm \mbox{for } t \; \mbox{in } (0,\infty)\\
\label{BCbis}
u(x,0)=u_0(x), \; \; \dfrac{\partial u}{\partial t}(x,0)=u_1(x) & & \hskip 1.2cm \mbox{for } x \; \mbox{in } (0,1) 
\label{ICbis}
\end{eqnarray}

\noindent Indeed A-R. Liu, C-H. Liu, J-Y. Fu, Y-L. Pi, Y-H. Huang and J-P. Zhang have realized experiments in which the value $k_{SMA}$ (equivalent stiffness of the damper) does not vanish. We start here with a simpler problem but we conjecture that the result will not change since the term we keep is the most important one for the stabilization. In this paper, $\beta=0$. The most important results are generalized to the case $\beta > 0$ in remarks (see Remarks \ref{rkeig} and \ref{rklarge}). \\ 
\\
Note that the case $a=1$, $b > 0$ and $\alpha=0$ is treated in \cite{AmDiZer}. They even consider $L^{\infty}(0;1)$ functions for $b$, which are positive and non-negative on an open subset of $(0;1)$. For a constant $b$, the energy is proved to decrease exponentially and the fastest decay rate is given by the supremum of the real part of the spectrum of the infinitesimal generator of the underlying semigroup. \\
\\
The case $a=1$, $b=0$ and $\alpha > 0$ is treated in \cite{AmTuc00} with two types of boundary conditions. They study the energy decay of a Bernoulli-Euler beam which is subject to a pointwise feedback force (given by the Dirac term). They show that both uniform and non uniform energy decay may occur depending on the boundary conditions. In the case of non uniform decay in the energy space (which is the case we are interested in), they give explicit polynomial decay estimates valid for regular initial data. Their method consists of deducing the decay estimates from observability inequalities for the associated undamped problem via sharp trace regularity results. \\
The same problem is studied by the same authors in \cite{AmTuc01} one year later. It is the second example of the applications given in Section 5. The estimates are unchanged for this example but the paper gives more general results: under a regularity assumption, the authors show that observability properties for the undamped problem imply decay estimates for the damped problem. 

\vspace{0.5cm}
\noindent Let $u$ be a regular solution of system (\ref{eq1})-(\ref{IC}). Its associated total energy is defined by

\begin{equation} \label{defE}
 E(t) = \dfrac{1}{2}\int_0^1(|u_t(x,t)|^2 + a |u_{xx}(x,t)|^2)dx.
\end{equation}

\noindent Then a classical computation using parts integration gives:
 
\begin{equation}
\dfrac{d}{dt}E(t) = -b \left(\int_0^1 |u_t(x,t)|^2 dx \right) -\alpha |u_t(\xi,t)|^2 \leq 0.
\end{equation}
Hence system (\ref{eq1})-(\ref{IC}) is dissipative in the sense that its associated energy is non increasing with respect to time.


\vspace{0.5cm}
\noindent A lot of studies have been performed by many authors in the stabilization of Euler-Bernoulli beams. Some of them can be found in the bibliography of the already cited papers : \cite{AmDiZer}, \cite{AmTuc00} and \cite{AmTuc01}. \\
The control of networks of Euler-Bernoulli or Timoshenko beams were studied in \cite{deko1}, \cite{caszu}, \cite{merreg1}, \cite{merreg2}, \cite{merreg3} and \cite{AkilChiGha} for example. Spectral methods are used like in this paper. \\
More recently, Euler-Bernoulli beams are coupled with wave equations with a Kelvin-Voigt damping in \cite{AkilAssaWehbe} for example. See also the references of this paper.

\noindent The main goal of this work is to obtain the energy decay rate of the damped Euler-Bernoulli beam subject to a pointwise feedback force representing a Shape Memory Alloy (SMA) cable described by (\ref{eq1}) to (\ref{IC}). \\
First we establish the well-posedness and strong stability of the problem. Then an explicit expression for the resolvent is given as well as the eigenvalues and eigenfunctions of the associated dissipative operator. We study the localization of the eigenvalues of the operator for small values of $\alpha > 0$. We prove that the eigenfunctions are quadratically close to those of the case $\alpha=0$. Thus they form a Riesz basis using a result of \cite{PTru}. \\
\\
At last, we give the explicit exponential decay estimate of the energy for regular initial data. The presence of a SMA cable ($\alpha > 0$) with a weak damping (small value for $\alpha$) may not restrain the vibrations of girder bridges. Our modelization does not take into account the stiffness term considered in \cite{LiuFu} ($\beta = 0$ here). The term which contributes most to the damping is $\alpha u_t$. That is why we start with this situation. See Remark \ref{rktaylor} which confirms this intuition. The stiffness term ($\beta \delta_{\xi} u$) even seems to diminish the decay rate of the energy. \\
\\
Note that the results were not obvious. One may think that adding a damping term ($\alpha > 0$ versus $\alpha = 0$) always leads to a better decay rate of the energy. This is not so clear. What we already know (cf. the conclusion of this paper) is that, if $\alpha = 0$, increasing the value of $b$ does not always lead to a higher decay rate of the energy, which is rather counterintuitive. \\  
\\
This paper does not confirm the results of \cite{LiuFu}. Nor does it refute them. The situation is more complicated. First they have an additional stiffness term but we conjecture that the results will be analogous. Indeed the eigenfunctions are unchanged and the characteristic equation is similar (cf. Remarks \ref{rkeig} and \ref{rklarge}). Moreover the large eigenvalues of the case $\alpha > 0$ and $\beta=0$ are proved to be close to those of the case $\alpha = 0$ and $\beta=0$ (see Proposition \ref{Largeeig}) and Remark \ref{rktaylor} proves that adding the stiffness term ($\beta >0$) leads to some eigenvalues with a real part larger than the supremum of the real part of the eigenvalues with $\beta = 0$. \\
All that is not necessarily in contradiction with the experiments of \cite{LiuFu}, since we have excluded some values for $\xi$ (cf. Theorem \ref{Theig}). As it is said in the conclusion of this paper, some more work is required to give a definite answer to the question in its title.   \rm 

\section{Well-posedness and strong stability} \label{sec2}

In this section we study the existence, uniqueness and strong stability of the solution of system (\ref{eq1})-(\ref{IC}).\\ 
\\
\noindent The energy space $\mathcal{H}$ is defined as follows

\be \label{defH}
\mathcal{H} = [H^2(0;1) \cap H_0^1(0;1)] \times L^2(0;1) 
\ee

\noindent with the inner product defined by

\beq
(U_1,U_2)_{\caH} = \displaystyle\int_0^1(a u_{1,xx} \overline{u_{2,xx}} + v_{1} \overline{v_{2}}) dx,\label{innprod}
\eeq

\noindent for all $U_1 = (u_1,v_1)$, $U_2 = (u_2,v_2)$ $\in$ $\caH$.\\
\\
\noindent Here again $a$ is a strictly positive constant (as in the introduction). \\

\brk
The norm $(U,U)_{\caH}^{\frac12}$ induced by (\ref{innprod}) is equivalent to the usual norm of $\caH$.
\erk

\noindent For shortness we denote by $\|\,.\,\|$   the  $L^2(\Omega)$-norm.

\noindent Now, we define a linear unbounded operator $\mathcal{A}:D(\mathcal{A})\rightarrow\caH$ by:

\begin{eqnarray} \label{defDA}
D(\mathcal{A}) = \bigg \{ U \in \mathcal{H} : (u,v) \in [H^4(0,\xi) \cap H^4(\xi,1) \cap H^2(0;1)] \times H^2(0;1),\, u(0)= v(0)=u(1)=v(1)=0 \\
\nonumber u_{xx}(0) = u_{xx}(1) = 0, \,  u_{xx}(\xi^-) = u_{xx}(\xi^+), \, u_x^{(3)}(\xi^+) - u_x^{(3)}(\xi^-) = -\dfrac{\alpha}{a} v(\xi) \bigg \} 
\end{eqnarray}

\begin{equation}
\mathcal{A}(u,v) = \big(v, - a u_{x}^{(4)} - b v - \alpha v(\xi) \delta_{\xi} \big),\ \ \ \  \forall\, U = (u,v) \in D(\mathcal{A}).\label{defA}
\end{equation}

\noindent Then we rewrite formally System (\ref{eq1})-(\ref{IC}) into the evolution equation

\begin{equation}
\left\{\begin{array}{ll}
U_t = \mathcal{A}U,\\
U(0) = U_0, \ \ \ \ U_0\in \caH
\label{abstr}
\end{array}\right.
\end{equation}
with  $U = (u,u_t)$.

\begin{proposition}
The operator $\caA$ is m-dissipative in the energy space $\caH$.
\end{proposition}

\proof On the first hand, the dissipativeness holds since we can check using integrations by parts: 

\begin{equation}
\Re (\caA U, U)_{\caH} = -b \left(\int_0^1 |u_t|^2 dx \right) -\alpha |u_t(\xi,t)|^2 \leq 0,  \forall \, U = (u,v) \in \caD(\caA). \label{AUU}
\end{equation} 

\noindent On the other end, the maximality is proved in the following way.

\noindent Let $f := (f_1;f_2) \in \caH$. We look for $U:=(u;v)  \in \mathcal{D}(\mathcal{A})$ solution of 

\begin{equation}\label{eqmaxmon}
- \mathcal{A} U =  f
\end{equation} 

or equivalently

\begin{equation}
\left\{\begin{array}{ll}
f_1 = - v \\
f_2 = - a u_{x}^{(4)} - b v - \alpha v(\xi) \delta_{\xi} 
\end{array} \right. \Leftrightarrow \left\{\begin{array}{ll}
v = -f_1 \\
a u_{x}^{(4)}  = f_2 + b f_1 - \alpha f_1(\xi) \delta_{\xi} 
\end{array} \right.
\end{equation}

\noindent Assume that such a solution $u$ exists, then multiplying the second identity by a function $\phi \in V:= H^2(0;1) \cap H_0^1(0;1)$, integrating in space and using integration by parts, it follows, since $U \in \mathcal{D}(\mathcal{A})$ 

\be \label{maxmoneq1}
\int_0^1 a u_x^{(2)} \phi_x^{(2)} dx = \int_0^1 (f_2 + b f_1) \phi dx.
\ee

\noindent This problem has a unique solution $u \in V = H^2(0;1) \cap H_0^1(0;1)$ by Lax-Milgram's lemma, because the left-hand side  of (\ref{maxmoneq1}) is coercive on $V$.

\noindent If we consider $\phi \in (\mathcal{D}(0,\xi)) \cap \mathcal{D}(\xi,1) \subset V$, then $u$ satisfies

$$a u_{x}^{(4)}  = f_2 + b f_1 - \alpha f_1(\xi) \delta_{\xi} \; \mbox{in the distributional sense}.$$

\noindent This directly implies that  $u \in H^4(0;\xi) \cap H^4(\xi;1)$ since $f_1$ and $f_2$ belong to $L^2(0;1)$. 

\noindent Now, integrating by parts lead to:
$$\begin{array}{ll} \label{energydecreasing}
\int_0^1 (f_2 + b f_1 - \alpha f_1(\xi) \delta_{\xi}) \phi dx = \int_0^{\xi} a u_x^{(4)} \phi dx + \int_{\xi}^{1} a u_x^{(4)} \phi dx +  
[a u_x^{(3)}(\xi^+) - a u_x^{(3)}(\xi^+)] \phi_x(\xi) \\
 + [a u_x^{(3)}(\xi^+) - a u_x^{(3)}(\xi^+)] \phi(\xi) + au_x^{(2)} \phi_x(1) - au_x^{(2)} \phi_x(0).
\end{array}$$

\noindent Consequently, by taking particular test functions  $\phi$ and using $v = -f_1$, we obtain $U \in \mathcal{D}(\mathcal{A})$ satisfying (\ref{eqmaxmon}), which finishes the proof of maximality.             \fin

\medskip

\noindent Using Lumer-Phillips Theorem (see \cite{P}, Theorem 1.4.3), the operator $\caA$ generates a $C_0$-semigroup of contractions $e^{t\caA}$ on
$\caH$. Then, we have the following result.

\bthe(Existence and uniqueness)\\
(1) If $U_0$ $\in$ $\caD(\caA)$, then   system (\ref {abstr}) has a unique strong solution
$$U  \in C^0(\R_+,\caD(\caA))\cap C^1(\R_+,\caH).$$
(2) If $U_0$ $\in$ $\caH$, then   system (\ref {abstr}) has a unique weak solution\\
$$U \in C^0(\R_+,\caH).$$
\ethe

\noindent Now, we have the following general strong stability result. 

\bthe \label{thstrong}(Strong stability)   \\
System (\ref{eq1})-(\ref{IC}) is strongly stable, i.e. for any solution $U$ of (\ref{abstr}) with initial data $U_0 \in\caH$, it holds 
$$\lim_{t\rightarrow\infty }E(t)=0,$$ where $E(t)$ is defined by (\ref{defE}). 
\ethe


\proof Since $\caA$ generates a contraction semigroup and its resolvent is compact in $\caH$, using Arendt-Batty Theorem(see \cite{AB}, p. 837), system (\ref{eq1})-(\ref{IC}) is strongly stable if and only if $\caA$ does not have purely imaginary eigenvalues. \\
\\
\noindent Assume that $\caA$ has a purely imaginary eigenvalue denoted by $i \mu$ with $\mu \in \R$. Then there exists $U = (u,v) \neq (0;0) \in \caD(\caA)$ such that $v= i \mu u$ and $\caA U = i \mu U$. Using (\ref{AUU}) leads to   

\begin{equation}
\Re \left[ i \mu (U, U)_{\caH} \right] = 0 = -b \left(\int_0^1 |v(x)|^2 dx \right) - \alpha |v(\xi)|^2. \label{AUU}
\end{equation} 

\noindent This implies $|v(\xi)| = 0$ and $\int_0^1 |v(x)|^2 dx=0$ (since $\alpha>0$ and $b > 0$). Thus $v \equiv 0$ almost everywhere on $(0;1)$ and since $v= i \mu u$, $u \equiv 0$ almost everywhere on $(0;1)$. Now $u$ and $v$ are continuous on $(0;1)$ and they vanish at $0$ and $1$ due to their belonging to $D(\caA)$. Thus $u \equiv v \equiv 0$. This contradicts the fact that $U \neq (0;0)$.   \fin


\section{Explicit expression for the resolvent} \label{sec3}
In this section we give an explicit expression of the resolvent $(\mu I -\caA)^{-1}$ and prove some useful estimates. In fact such estimates are useful since later on, we will use a result of \cite{Bor} (Theorem 2.4) which involves the norm operator of $(\mu I -\caA)^{-1}$  with  $\mu \in \C.$ 

\noindent Let $F=(u_1,v_1) \in \caH,$ we look for a solution $U =(u,v) \in \da$  of 
\be \label{r1}( \mu I - \caA) U= F , \; \mu \in \C.
\ee 

The explicit expression for the resolvent we give in next Proposition \ref{propresol} involves the restriction on $[0;1]$ of the classical convolution product of two functions on $\R$. Let us recall the definition and two useful properties we established in \cite{EECT18}.   

\blem \label{convol} (A technical lemma) \\ 
Let  $\psi \in C^{\infty}([0,\infty[)$ and $f$ in $L^2(0; 1)$ be two functions and define their convolution product $\psi \star f$ on $[0,1]$ by :

\be \label{convolution1} 
(\psi \star f)(x) = \int_0^x \psi(x-s) f(s) ds, \forall x \in [0;1].
\ee

Then the following two properties hold:
\begin{enumerate}
\item $(\psi \star f) \in H^1(0; 1)$ and its derivative is:
\be \label{prop1}
(\psi \star f)'(x) = \int_0^x \psi'(x-s) f(s) ds + \psi(0)f(x), \forall x \in [0; 1].
\ee

\item If $\psi(0) = 0$ is also assumed, then $(\psi \star f) \in H^2(0; 1)$ and its second derivative is:

\be \label{prop2}
(\psi \star f)''(x) = \int_0^x \psi''(x-s) f(s) ds + \psi'(0)f(x), \forall x \in [0; 1].
\ee

\end{enumerate} 
\elem 
\proof 
\begin{enumerate}
\item The functions $\psi$ extended by $0$ on $(-\infty,0)$ and $f$  extended by $0$ on $\R$  outside $[0,1]$ are still called $\psi$ and $f$ respectively. Then the convolution product defined by (\ref{convolution1}) is extended by the classical convolution product on $\R$ i.e by 

\be \label{convolution2} 
(\psi \star f)(x)=\int _R  \psi(x-s) f(s) ds, \forall x \in \R. \ee 

It is well known that $(\psi \star f)' =(\psi ')_{\mbox{dist}} \star f$ where $(\psi ')_{\mbox{dist}}$ is the derivative of $\psi$ in the distributional sense. 
Due to the property of $\psi$ and its extension on $\R$ we have 
$$(\psi ')_{\mbox{dist}}=\psi'+\psi(0)\delta_0,$$ where $\delta_0$ is the Dirac distribution at $x=0.$ The property (\ref{prop1}) follows from this remark. 
 
\item (\ref{prop2}) is a consequence of (\ref{prop1}). \fin
\end{enumerate}

\noindent Note that weaker assumptions could be made on $\psi$ for this lemma ($\psi \in C^2([0,1])$ is sufficient).

\bprop \label{propresol} (Explicit expression for the resolvent of the operator $\caA$) \\
Let $a$, $b$ and $\alpha$ be strictly positive real numbers, $\xi$ a real number in $(0;1)$ and $\mu$ a complex number. \\
Let the spaces $\caH$ and $D(\caA)$ be defined by $(\ref{defH})$ and $(\ref{defDA})$. Let $F=(u_1,v_1) \in \caH$. \\
Denote by $\lambda$ the complex number, such that $\exists k \in \Z, \arg(\lambda) \in [-\pi/4 + 2 k \pi;\pi/4 + 2 k \pi)$, satisfying: 

\be \label{deflambda}
\lambda^4 = -\dfrac{b \mu + \mu^2}{a}.
\ee 

\noindent Denote by $H(x,0)$ the Heaviside step function defined by 

\be \label{H0}
H(x,0) : = \left\{\begin{array}{ll}
0, \forall x < 0 \\
1, \forall x \geq 0.  
\end{array} \right.
\ee

\noindent Define the expressions:

\be \label{f10}
f_1^0:= (\mu + b)u_1 + v_1
\ee

\be \label{defDetalpha}
Det^{\alpha}(\lambda):= 4 \lambda^2 \left \{ [-2 \lambda^3 \sinh(\lambda) - \dfrac{\alpha \mu}{a} \sinh(\lambda \xi) \sinh(\lambda (\xi - 1))] \sin(\lambda) + \dfrac{\alpha \mu}{a} \sin(\lambda \xi) \sinh(\lambda) \sin(\lambda (\xi - 1))   \right \}
\ee

\be \label{defDet0}
Det^0(\lambda):= - 8 \lambda^5 \sinh(\lambda) \sin(\lambda)
\ee

\be \label{u0}
u_0(\lambda,x):= \dfrac{1}{2 a \lambda^3} \left[ \sin(\lambda x) - \sinh(\lambda x) \right] H(x,0), \forall x \in (0;1)
\ee

\be \label{ABC}
\left\{\begin{array}{ll}
A(\lambda, \alpha) := \lambda^2 (u_0(\lambda, \cdot) \star f_1^0)(1) + 2 \lambda^2 \alpha u_1(\xi) u_0(\lambda, 1-\xi)  \\ 
B(\lambda, \alpha) := (u_0(\lambda, \cdot) \star f_1^0)_{xx}(1) + 2 \alpha u_1(\xi) u_{0,xx}(\lambda, 1-\xi) \\
C(\lambda,\xi) := \lambda^2 (u_0(\lambda, \cdot) \star f_1^0)(\xi) 
\end{array} \right.
\ee

\noindent For simplicity, the dependency of $A$, $B$ and $C$ on $\lambda$ and $\xi$ is omitted in the following.

\be \label{resolvent0}
\begin{array}{ll}
R(\mu, \caA_0)F : = (u_0 \star f_1^0)(x) + \dfrac{a}{Det^0(\lambda)} \bigg[& 4 \lambda^4 \left \{(A-B)\sinh(\lambda) + (A+B)\sin(\lambda) \right \} u_0(x,\lambda) \\
& + 4 \lambda^6 \left \{(A+B)\sin(\lambda) - (A-B)\sinh(\lambda) \right \} u_{0,xx}(x,\lambda)  
\bigg]
\end{array}
\ee

\be \label{defDelta}
\Delta(\lambda,\xi):= \dfrac{4 \lambda^2 \mu}{a} \left( \sin(\lambda \xi) \sin(\lambda(\xi-1)) \sinh(\lambda) - \sinh(\lambda \xi) \sinh(\lambda(\xi-1)) \sin(\lambda) \right)
\ee

\be \label{defDelta3}
\begin{array}{ll}
\Delta_3(\lambda,\xi):= \dfrac{4 \lambda \mu}{a} \bigg \{ 
& \sin(\lambda (\xi - 1)) \left( \dfrac{A+B}{2} [\sinh(\lambda \xi) - \sin(\lambda \xi)] - C \sinh(\lambda) \right) \\
& + \sinh(\lambda(\xi-1)) \left( \dfrac{A-B}{2} [\sinh(\lambda \xi) - \sin(\lambda \xi)] + C \sin(\lambda) \right)
\bigg \}
\end{array}
\ee

\be \label{defDelta1}
\begin{array}{ll}
\Delta_1(\lambda,\xi):= \dfrac{4 \lambda^3 \mu}{a} \bigg \{ 
& \sin(\lambda (\xi - 1)) \left(  - \dfrac{A+B}{2} [\sinh(\lambda \xi) + \sin(\lambda \xi)] + C \sinh(\lambda) \right) \\
& + \sinh(\lambda(\xi-1)) \left( - \dfrac{A-B}{2} [\sinh(\lambda \xi) + \sin(\lambda \xi)] + C \sin(\lambda) \right)
\bigg \}
\end{array}
\ee

\be \label{defDelta0}
\Delta_0(\lambda,\xi):= - 4 \lambda^3  \bigg \{ 
2 C \sin(\lambda) \sinh(\lambda) - (A + B) \sin(\lambda) \sinh(\lambda \xi) - (A - B) \sin(\lambda \xi) \sinh(\lambda) \bigg \} 
\ee

\noindent Then the solution $U=(u,v) \in \da$ of $( \mu I - \caA) U= F$, with $F=(u_1,v_1) \in \caH,$ is given by:  

\be \label{resolvent}
\begin{array}{ll}
u(x)= & \dfrac{1}{Det^{\alpha}(\lambda)} \bigg(Det^0(\lambda) \cdot R(\mu, \caA_0)F + \alpha \bigg[ \Delta(\lambda,\xi) (u_0 \star f_1^0)(x) + a \Delta_3(\lambda,\xi)u_0(\lambda,x) + a \Delta_1(\lambda,\xi)u_{0,xx}(\lambda,x) \\
& - \mu \Delta_0(\lambda, \xi)u_0(\lambda,x - \xi) \bigg] \bigg)
\end{array}
\ee

\eprop 

\proof The computation of an explicit expression for the resolvent of $\caA$ requires the search for a solution $u$ in $H^4(0;\xi) \cap H^2(0;1)$ (resp. in $H^4(\xi;1) \cap H^2(0;1)$) of the following problem on $(0;\xi)$ (resp. on $(\xi;1)$):

\be \label{pbresol}
\left \{ \begin{array}{llllll}
a u_x^{(4)} + (b \mu + \mu^2) u + \mu \alpha \delta_{\xi} u = f_1^{\alpha} \\
v=\mu u - u_1 \\
u(0)=u(1)=0 \\
u_{xx}(0) = u_{xx}(1)=0 \\
u_{xx}(\xi^-) = u_{xx}(\xi^+) \\
u_{x}^{(3)}(\xi^+) - u_{x}^{(3)}(\xi^-)=g_1
\end{array} \right.
\ee

\noindent with $f_1^{\alpha}:=(\mu + b)u_1 + v_1 + \alpha \delta_{\xi} u_1$ and $g_1:= \dfrac{\alpha}{a} (- \mu u(\xi) + u_1(\xi))$. Note that $f_1^0$ defined by (\ref{f10}) is $f_1^{\alpha}$ with $\alpha=0$.

\noindent First we still denote by $u$ the function defined by

\be 
u(x) : = \left\{\begin{array}{ll}
0, \forall x < 0 \; \mbox{and} \; x > \xi \\
u(x), \forall x \in (0;\xi).  
\end{array} \right.
\ee

\noindent where $u$ is a solution of (\ref{pbresol}) on $(0;\xi)$. 

\noindent Denote by $\caL$ the classical Laplace transform. Using four successive integrations by parts and the properties of $u$ at $0$, $1$, $\xi^-$ and $\xi^+$ leads to:

\be
[\caL u_x^{(4)}](p) = p^4 [\caL u](p) + [u_x^{(3)}(\xi^-) + p u_x^{(2)}(\xi^-)+ p^2 u_x(\xi^-)+p^3 u(\xi^-)]e^{-p \xi} - [u_x^{(3)}(0) + p^2 u_x(0)] 
\ee

\noindent Applying the Laplace transform to Problem (\ref{pbresol}) on $(0;\xi)$ gives:

\be
[a p^4 + b \mu + \mu^2] [\caL u](p) + \alpha \mu u(\xi) e^{-p \xi} + a [u_x^{(3)}(\xi^-) + p u_x^{(2)}(\xi^-)+ p^2 u_x(\xi)+p^3 u(\xi)]e^{-p \xi}  - a [u_x^{(3)}(0) + p^2 u_x(0)] = \caL (f_1^{\alpha} \cdot \chi_{(0;\xi)}) 
\ee

\noindent We proceed similarly on $(\xi;1)$ and sum both functions to get:

\be \label{Laplaceu}
\begin{array}{lll}
[\caL u](p) = \dfrac{1}{ap^4 + b \mu + \mu^2} \bigg[ [\caL f_1^{\alpha}](p) + a u_x^{(3)}(0) + (u_x^{(3)}(\xi^+) - u_x^{(3)}(\xi^-)) a e^{-p \xi} - a u_x^{(3)}(1) e^{-p} \\
+ a p^2 [u_x(0) - u_x(1) e^{-p}] - \alpha \mu u_1(\xi) e^{-p \xi} \bigg] \\
= \dfrac{1}{ap^4 + b \mu + \mu^2} \bigg[ [\caL f_1^{\alpha}](p) + a u_x^{(3)}(0) + \alpha u(\xi) e^{-p \xi} - a u_x^{(3)}(1) e^{-p} + a p^2 [u_x(0) - u_x(1) e^{-p}] \bigg] 
\end{array}
\ee

\noindent The function $u_0$ defined by (\ref{u0}) satisfies: $\left[ \caL \left[ \dfrac{1}{2 a \lambda^3} \left[ \sin(\lambda \cdot) - \sinh(\lambda \cdot) \right] \right] \right](p) = \dfrac{1}{ap^4 + b \mu + \mu^2}$. Thus 

\be \label{u}
u(x) = [u_0(\lambda, \cdot) \star f_1^{\alpha}](x) + a u_x^{(3)}(0) \cdot u_0(\lambda,x) + a u_x(0) \cdot u_{0,xx}(\lambda,x) + [\alpha u_1(\xi) - 2 \alpha \mu u(\xi)] \cdot u_0(\lambda,x - \xi)
\ee

\noindent An analogous problem is treated in Section $2.1$ of \cite{Maj} for example. Taking into account $u(1)= u_{xx}(1)=0$ gives the following two equations:

\be \label{sys12}
\left \{ \begin{array}{llll}
\lambda^2 [\sinh(\lambda) + \sin(\lambda)] u_x^{(3)}(0) + [\sinh(\lambda) - \sin(\lambda)] u_x(0) - \dfrac{\alpha \mu}{a} [\sinh(\lambda(1-\xi)) - \sin(\lambda (1-\xi))] u(\xi) \\
= - 2 \lambda^3 (u_0 \star f_1^0)(1) - 4 \lambda^3 \alpha u_1(\xi) u_0(1-\xi) \\
\lambda^2 [\sinh(\lambda) - \sin(\lambda)] u_x^{(3)}(0) + [\sinh(\lambda) + \sin(\lambda)] u_x(0) - \dfrac{\alpha \mu}{a} [\sinh(\lambda(1-\xi)) + \sin(\lambda (1-\xi))] u(\xi) \\
= - 2 \lambda (u_0 \star f_1^0)^{(2)}(1) - 4 \lambda \alpha u_1(\xi) u_{0,xx}(1-\xi)  
\end{array} \right.
\ee

\noindent (\ref{prop2}) has been used here since $u_0(0)=u_{0,x}(0)=0$. The last step is to evaluate (\ref{u}) at $x=\xi$: 

\be \label{sys3}
[\sinh(\lambda \xi) - \sin(\lambda \xi)] u_x^{(3)}(0) + \lambda^2 [\sinh(\lambda \xi) + \sin(\lambda \xi)] u_x(0) - 2 \lambda^3 u(\xi) = - 2 \lambda^3 (u_0 \star f_1^0)(\xi).  
\ee

\noindent Now (\ref{sys12}) and (\ref{sys3}) form a $3 \times 3$ system where the unknowns are $u_x^{(3)}(0)$, $u_x(0)$ and $u(\xi)$. Its determinant is $Det^{\alpha}(\lambda)$. Solving this system gives (\ref{resolvent}) after long calculations. \fin

\brk The expression $R(\lambda, \caA_0)$ represents the resolvent of the operator $\caA_0$ which is $\caA$ with $\alpha=0$. Thus the expression for the solution $u$ given by (\ref{resolvent}) is $R(\lambda, \caA_0)F$ if $\alpha=0$.
\erk


\section{Eigenvalues and eigenfunctions of the operator $\caA$} \label{sec4}

The eigenvalues and eigenfunctions of the operator $\caA$ defined by (\ref{defDA}) and (\ref{defA}) in Section \ref{sec2} are explicitly calculated. For the sake of completeness, the case $\alpha = 0$ which corresponds to the operator $\caA_{a_0}$ of \cite{AmDiZer} is recalled here. Note the dependency of $\caA$ on $a$, $b$, $\alpha$ and $\xi$. \\
\\
First a technical Lemma will be useful in the proof of Theorem \ref{Theig} to handle with particular values for $\xi$.

\blem (Is $\mu = - b$ an eigenvalue of $\caA$?) \label{lemb}\\
Assume that $a$, $b$, $\alpha$ are strictly positive constants, $\xi \in (0;1)$. The following problem on $H^4(0;\xi)\cap H^2(0;1)$ (resp. on $H^4(\xi;1) \cap H^2(0;1)$):

\be \label{beig}
\left \{ \begin{array}{lllll}
a u_x^{(4)} - b \alpha \delta_{\xi} u = 0 \\
u(0)=u(1)=0 \\
u_{xx}(0) = u_{xx}(1)=0 \\
u_{xx}(\xi^-) = u_{xx}(\xi^+) \\
u_{x}^{(3)}(\xi^+) - u_{x}^{(3)}(\xi^-)= \dfrac{\alpha}{a} b u(\xi)
\end{array} \right.
\ee

\noindent has no other solution than $u \equiv 0$ if the polynomial $P(x)= 1 - \dfrac{2 \alpha b}{3} x^2 (1-x)^2$ has no solution in $(0;1)$. There are three cases:

\begin{itemize}
\item Case $\alpha b > 6$: $P$ has exactly two roots in $(0;1)$ called $\xi_1$ and $\xi_2$ in the following.
\item Case $\alpha b = 6$: $P$ has exactly one root in $(0;1)$ called $\xi_0$ in the following. 
\item Case $\alpha b < 6$: the polynomial $P$ does not vanish on $(0;1)$.
\end{itemize}

\elem

\proof The proof starts like that of Proposition \ref{propresol}. Denote by $u_0$ the function defined by $u_0(x)= \dfrac{x^3}{6a} \cdot H(x,0), \forall x \in (0;1)$. It holds $[\caL \left( \dfrac{x^3}{6a} \right)](x \mapsto p)= \dfrac{1}{a p^4}, \forall p \in (0;+ \infty)$. Thus     

\be \label{newu}
u(x) = a u_x^{(3)}(0) \cdot u_0(\lambda,x) + a u_x(0) \cdot u_{0,xx}(\lambda,x) + \dfrac{\alpha}{a} b u(\xi) \cdot u_0(\lambda,x - \xi).
\ee

\noindent Taking into account $u(1)= u_{xx}(1)=0$ gives the following two equations:

\be \label{newsys12}
\left \{ \begin{array}{ll}
\dfrac{1}{6} u_x^{(3)}(0) + u_x(0) + \dfrac{\alpha}{a} b u(\xi) (1 - \xi)^3 = 0 \\
u_x^{(3)}(0) + 2 \alpha b u(\xi) (1 - \xi) = 0
\end{array} \right.
\ee

\noindent The last step is to evaluate (\ref{newu}) at $x=\xi$: 

\be \label{newsys3}
\dfrac{\xi^3}{6} u_x^{(3)}(0) + \xi u_x(0) - u(\xi) = 0. 
\ee

\noindent Now (\ref{newsys12}) and (\ref{newsys3}) form a $3 \times 3$ system where the unknowns are $u_x^{(3)}(0)$, $u_x(0)$ and $u(\xi)$. An obvious solution is $(0;0;0)$. Its determinant is $P(\xi)$. \\
The system has no other solution than $(0;0;0)$ if and only if $P(\xi)$ does not vanish. \\
The roots of $P$ follow from a classical study of the variations of the function $P$.
\fin 


\bthe (Eigenvalues and eigenfunctions of the operator $\caA$) \label{Theig}\\
Let $\caA$ be the operator defined by (\ref{defDA}) and (\ref{defA}) in Section \ref{sec2}. Denote by $\mu_n$, $n \in \Z$ the eigenvalues of $\caA$ and by $\lambda$ the complex numbers satisfying $\lambda^4 = - \dfrac{b \mu + \mu^2}{a}$. One of these $\lambda$'s is such that $\exists k \in \Z, \arg(\lambda) \in [-\pi/4 + 2 k \pi;\pi/4 + 2 k \pi)$. The others are $- \lambda$ and $\pm i \lambda$ and they are still denoted by $\lambda$.

\begin{enumerate}
\item Case $\alpha=0$. \\
If $b \in \R^{+*}- \{ 2 \sqrt{a} n^2 \pi^2, n \in \Z^* \}$, then $\exists n_0$ such that $2 \sqrt{a} n_0^2 \pi^2 < b < 2 \sqrt{a} (n_0 + 1)^2 \pi^2$ and 

\be \mu_n^{\pm} = \left \{ \begin{array}{ll}
\dfrac{1}{2} \left( -b \pm \sqrt{b^2 - 4 a n^4 \pi^4} \right), n = 1,2, \ldots, n_0 \\
\dfrac{1}{2} \left( -b \pm i \sqrt{4 a n^4 \pi^4 - b^2} \right), n > n_0. 
\end{array}
\right.
\ee 
\noindent The algebraic multiplicity of these eigenvalues is $1$.\\
\\
If $\exists n_0$ such that $b = 2 \sqrt{a} n_0^2 \pi^2$ then $\mu_{n_0}= - \dfrac{b}{2}$ and the algebraic multiplicity of this eigenvalue may not exceed $2$.\\
\\
\noindent In both cases, the associated eigenfunction is defined on $(0;1)$ by 

\be \label{eigfunc0}
\Phi_n^{\pm}(x)= \sin(n \pi x)(1, \mu_n^{\pm}), n \geq 1. 
\ee

\item Case $\alpha \neq 0$. \\

   \begin{enumerate}
      \item Case $\alpha b < 6$: $\mu = -b$ is not an eigenvalue of $\caA$. \\
\noindent If $\xi \notin \Q$, $\mu$ is an eigenvalue of $\caA$ if and only if the corresponding $\lambda$'s satisfy $\lambda \notin \{ k \pi, k \in \Z \}$, $\lambda \notin \{ i k \pi, k \in \Z \}$ and the characteristic equation:

\be \label{eqcar}
2 (\mu + b) \sinh(\lambda) \sin(\lambda) + \alpha \lambda \left[
\sin(\lambda) \sinh(\lambda \xi) \sinh(\lambda(1 - \xi))- \sinh(\lambda) \sin(\lambda \xi) \sin(\lambda(1 - \xi)) 
\right] = 0
\ee

\noindent and the associated eigenfunction is $\Phi(\mu, \cdot)= (\phi(\lambda, \cdot)(1,\mu)$ where $\phi(\lambda, \cdot)$ is defined on $(0;1)$ by 

\be \label{eigfunc}
\begin{split}
\phi(\lambda, x):= & \dfrac{1}{|\lambda|^2} \cdot e^{(\xi - 2) |\lambda|} \bigg \{ \sin(\lambda) \sinh(\lambda) [\sin(\lambda(x - \xi)) - \sinh(\lambda(x - \xi))] \cdot H(x, \xi) \\
&  + \sinh(\lambda) \sin(\lambda(1-\xi)) \sin(\lambda x) - \sin(\lambda) \sinh(\lambda(1-\xi)) \sinh(\lambda x) \bigg \} 
\end{split}
\ee

\noindent with $H(x,\xi)$ the Heaviside step function at $\xi$. \\
\\
Denote by $r = r(a,b)$ the real positive number such that $r^4= \dfrac{b^2}{4a}$. Then, if  

\be \label{partcase}
\alpha = - \dfrac{4 a r^3 \sinh(r) \sin(r)}{b [\sinh(r) \sin(r \xi) \sin(r(\xi - 1)) - \sin(r) \sinh(r \xi) \sinh(r(\xi - 1))]} 
\ee
then $\mu=-b/2$ is an eigenvalue of $\caA$ with algebraic multiplicity at least $2$ and its geometric multiplicity is $1$. \\
\\
If $\alpha$ takes any other value, the algebraic multiplicity of $\mu$ as an eigenvalue of $\caA$ is at least $1$ and its geometric multiplicity is $1$.
\item Case $\alpha b = 6$: there exists $\xi_0$ in $(0;1)$ such that, if $\xi=\xi_0$, $\mu = -b$ is an eigenvalue of $\caA$ ($=\caA(\xi)$). If $\xi \neq \xi_0$, $\mu = -b$ is not an eigenvalue of $\caA$ ($=\caA(\xi)$). In both cases, if $\xi \notin \Q$, the other eigenvalues and the associated eigenfunctions are given by the same expressions as in the preceding case.
\item Case $\alpha b > 6$: there exist $\xi_1$ and $\xi_2$ in $(0;1)$ such that, if $\xi=\xi_1$ or $\xi=\xi_2$, $\mu = -b$ is an eigenvalue of $\caA$ ($=\caA(\xi)$). If $\xi \notin \{ \xi_1, \xi_2 \}$, $\mu = -b$ is not an eigenvalue of $\caA$ ($=\caA(\xi)$). In both cases, if $\xi \notin \Q$, the other eigenvalues and the associated eigenfunctions are given by the same expressions as in the first case. 
   
\end{enumerate}

\end{enumerate}

\ethe  

\proof First of all, the operator $\caA$ has already been proved to have no imaginary eigenvalues (cf. the proof of Theorem \ref{thstrong}). In particular, $0$ is not an eigenvalue. \\
\\
\noindent The computation of the eigenelements of $\caA$ requires the search for a solution $u$ in $H^4(0;\xi) \cap H^2(0;1)$ (resp. in $H^4(\xi;1) \cap H^2(0;1)$) of the following problem on $(0;\xi)$ (resp. on $(\xi;1)$):

\be \label{eigpb}
\left \{ \begin{array}{llllll}
a u_x^{(4)} + (b \mu + \mu^2) u + \mu \alpha \delta_{\xi} u = 0 \\
u(0)=u(1)=0 \\
u_{xx}(0) = u_{xx}(1)=0 \\
u_{xx}(\xi^-) = u_{xx}(\xi^+) \\
u_{x}^{(3)}(\xi^+) - u_{x}^{(3)}(\xi^-)= - \dfrac{\alpha}{a} \mu u(\xi)
\end{array} \right.
\ee

\vspace{0.5cm}
\bf Case $\alpha = 0$: \rm the proof is that of Lemma 2.2 of \cite{AmDiZer}, mutatis mutandis.
\vspace{0.5cm}

\bf Case $\alpha \neq 0$ and $\alpha b < 6$: \rm since neither $\mu =0$ nor $\mu = -b$ is an eigenvalue of $\caA$ (cf. Lemma \ref{lemb}), $\lambda$ cannot vanish and a particular solution of Problem (\ref{eigpb}) is:

\be \label{up}
u_p(\lambda,x):= - \dfrac{\alpha \mu}{2 a \lambda^3} u(\xi) \left[ \sinh(\lambda (x-\xi)) - \sin(\lambda (x-\xi)) \right] H(x,\xi), \forall x \in (0;1)
\ee

\noindent where $H(x,\xi)$ is for the Heaviside step function defined by 

\be \label{H0}
H(x,\xi) : = \left\{\begin{array}{ll}
0, \forall x < \xi \\
1, \forall x \geq \xi.  
\end{array} \right.
\ee

\noindent This result is not new: it can be found in \cite{Maj} for example.\\
\\
\noindent The general solution of this problem (i.e. with $u_{x}^{(3)}(\xi^+) - u_{x}^{(3)}(\xi^-)=0$) can be written as:

\be \label{ug}
u_g(\lambda,x):= \left \{ \begin{array}{ll}
 P_1 \cosh(\lambda x) + Q_1 \sinh(\lambda x) + R_1 \cos(\lambda x) + S_1 \sin(\lambda x), \forall x \in (0;\xi] \\
 P_2 \cosh(\lambda (x - \xi)) + Q_2 \sinh(\lambda (x - \xi)) + R_2 \cos(\lambda (x - \xi)) + S_2 \sin(\lambda (x - \xi)), \forall x \in [\xi;1) 
\end{array} \right.
\ee

\noindent Using the boundary conditions for $u$ and $u_{xx}$ as well as the continuity of $u_{xx}$ and $u_{x}^{(3)}$ at $\xi$, it holds:

\be \label{pqrs}
\left \{ \begin{array}{lll}
P_1 = R_1 = 0 \\
P_2 = Q_1 \sinh(\lambda \xi), Q_2 = Q_1 \cosh(\lambda \xi), R_2 = S_1 \sin(\lambda \xi) \; \mbox{and} \; S_2 = S_1 \cos(\lambda \xi) \\
\sinh(\lambda) Q_1 = - \dfrac{\alpha \mu}{2 a \lambda^3} u(\xi) \sinh(\lambda (\xi-1)) \; \mbox{and} \; \sin(\lambda) S_1 =  \dfrac{\alpha \mu}{2 a \lambda^3} u(\xi) \sin(\lambda (\xi-1)) 
\end{array} \right. 
\ee

\vspace{0.5cm}

If $\xi \notin \Q$, $\lambda \notin \{ k \pi, k \in \Z^* \}$ and $\lambda \notin \{ i k \pi, k \in \Z^* \}$. Indeed, if $\exists k \in \Z^*$, such that $\lambda = k \pi$, then $\sin(\lambda)=0$ but $\sinh(\lambda) \neq 0$ and the last identity of (\ref{pqrs}) is:

\be 
\dfrac{\alpha \mu}{2 a \lambda^3} u(\xi) \sin(k \pi (\xi-1)) = 0  
\ee

Since $\alpha \neq 0$, $\mu \neq 0$ and $\sin(k \pi (\xi - 1)) \neq 0$ (the expression only vanishes if $\xi \in \Q$), it holds $u(\xi)=0$. Now, since $\sinh(\lambda) \neq 0$, it implies $Q_1=0$ and consequently, $P_2=Q_2=0$. The vanishing of $u(\xi)$ also implies that of $u_p$. Thus

\be 
u(\lambda, x) = u_g(\lambda,x):= \left \{ \begin{array}{ll}
S_1 \sin(n \pi x), \forall x \in (0;\xi] \\
S_1 \sin(n \pi \xi) \cos(n \pi (x - \xi)) + S_1 \cos(n \pi \xi) \sin(n \pi (x - \xi)) = S_1 \sin(n \pi x), \forall x \in [\xi;1). 
\end{array} \right.
\ee

Now, $u(\xi)=0=S_1 \sin(n \pi \xi)$ and $S_1 \neq 0$ (otherwise $u \equiv 0$), then $\sin(n \pi \xi) = 0$, which contradicts the fact that $\xi \notin \Q$. \\
\\
Analogously, $\lambda \notin \{ i k \pi, k \in \Z^* \}$ if $\xi \notin \Q$. \\
\\
Now, since $\lambda \notin \{ k \pi, k \in \Z^* \}$ and $\lambda \notin \{ i k \pi, k \in \Z^* \}$, (\ref{pqrs}) becomes:  

\be \label{pqrs3}
\left \{ \begin{array}{lll}
P_1 = R_1 = 0 \\
Q_1 = - \dfrac{\alpha \mu}{2 a \lambda^3} u(\xi) \dfrac{\sinh(\lambda (\xi-1))}{\sinh(\lambda)}  \; \mbox{and} \; S_1 = \dfrac{\alpha \mu}{2 a \lambda^3} u(\xi) \dfrac{\sin(\lambda (\xi-1))}{\sin(\lambda)} \\
P_2 = Q_1 \sinh(\lambda \xi), Q_2 = Q_1 \cosh(\lambda \xi), R_2 = S_1 \sin(\lambda \xi) \; \mbox{and} \; S_2 = S_1 \cos(\lambda \xi) 
\end{array} \right. 
\ee 

\noindent Combining that with classical trigonometric formulae leads to: 

\be 
u_g(\lambda,x):= Q_1 \sinh(\lambda x) + S_1 \sin(\lambda x), \forall x \in (0;1).
\ee

\noindent In particular $u(\xi)= u_p(\lambda, \xi) + u_g(\lambda, \xi) = u_g(\lambda, \xi)$ i.e.

\be 
u(\xi)= - \dfrac{\alpha \mu}{2 a \lambda^3} u(\xi) \left( \dfrac{\sinh(\lambda (\xi-1))}{\sinh(\lambda)} \sinh(\lambda x) - \dfrac{\sin(\lambda (\xi-1))}{\sin(\lambda)} \sin(\lambda x) \right).
\ee

\noindent This implies the following equation, since $u(\xi) \neq 0$:

\be \label{eqcar1}
2 a \lambda^3 \sinh(\lambda) \sin(\lambda) - \alpha \mu \left[
\sin(\lambda) \sinh(\lambda \xi) \sinh(\lambda(1 - \xi))- \sinh(\lambda) \sin(\lambda \xi) \sin(\lambda(1 - \xi)) 
\right] = 0
\ee

\noindent Multiplying both sides by $\lambda \neq 0$, replacing $(a \lambda^4)$ by $(- b \mu - \mu^2)$ and dividing by $\mu \neq 0$ leads to the characteristic equation (\ref{eqcar}). \\
\\
\noindent At last, $u(\xi)$ can be rewritten as:

\be 
u(\xi)= u_g(\xi) = - \dfrac{\sinh(\lambda (\xi-1))}{\sinh(\lambda)} \sinh(\lambda x) + \dfrac{\sin(\lambda (\xi-1))}{\sin(\lambda)} \sin(\lambda x) = \dfrac{2 a \lambda^3}{\alpha \mu} \; \mbox{due to the characteristic equation.}
\ee

\noindent Hence the expression for the eigenfunction (\ref{eigfunc}). Note that the factor $\dfrac{1}{|\lambda|^2} e^{(\xi - 2) |\lambda|}$ is aimed at making the function bounded with $\lambda$ in $\caH$. \\
\\
\noindent The last part of the proof concerns the multiplicity of the eigenvalues. Since the characteristic equation (\ref{eqcar}) is invariant under the transform $\lambda \mapsto i \lambda$, to each root $\lambda$ such that $\lambda \in \{ z \in \C - \{ k \pi, k \in \Z \},  \exists l \in \Z, \arg(z) \in [-\pi/4 + 2 l \pi;\pi/4 + 2 l \pi)$ correspond three other roots which are $\pm i \lambda$ and $- \lambda$. Up to a multiplicative constant, the expression for the eigenfunction is invariant under the transform $\lambda \mapsto i \lambda$, thus the geometric multiplicity of $\mu$ as an eigenvalue is always $1$. \\
\\
\noindent Now (\ref{deflambda}) is equivalent to 

\be \label{mu}
\mu(\lambda) = \dfrac{-b + \delta(\lambda)}{2}
\ee

\noindent where $\pm \delta(\lambda)$ are the (complex) square roots of $(b^2 - 4 a \lambda^4)$. \\
\\
The value $\mu = -b/2$ is an eigenvalue of the operator $\caA$ if and only if $\delta(\lambda)=0$ which is the case for the value $\alpha$ defined by (\ref{partcase}). In that case, the algebraic multiplicity is at least $2$. \\
If $\alpha$ takes any other value, it is at least $1$.\fin
   

\brk The characteristic equation (\ref{eqcar}) (found for $\alpha \neq 0$ and $\xi \notin \Q$) becomes $\sinh(\lambda) \sin(\lambda)=0$, if $\alpha$ is replaced by $0$ in (\ref{eqcar}). The solutions of this equation are $\lambda = n \pi$, $n \in \Z$ and $\lambda = i n \pi$, $n \in \Z$. Now, if $\alpha$ tends to $0$, $\alpha b < 6$ so $\mu = - b$ is not an eigenvalue of $\caA$ and since $0$ is not an eigenvalue of $\caA$, $\lambda \neq 0$. Thus the corresponding eigenvalues are those given for Case $\alpha=0$ in Theorem \ref{Theig}, which is coherent.
\erk

\brk (Case $\beta>0$) \label{rkeig} The eigenfunctions of Problem (\ref{eq1bis})-(\ref{ICbis}) with the additional term $\beta \delta_{\xi}$ are unchanged (cf. (\ref{eigfunc})) and the characteristic equation (\ref{eqcar}) becomes:   

\be 
2 a \lambda^3 \sinh(\lambda) \sin(\lambda) - (\alpha \mu + \beta) \left[
\sin(\lambda) \sinh(\lambda \xi) \sinh(\lambda(1 - \xi))- \sinh(\lambda) \sin(\lambda \xi) \sin(\lambda(1 - \xi)) 
\right] = 0.
\ee

\noindent Multiplying both sides by $\lambda \neq 0$, replacing $(a \lambda^4)$ by $(- b \mu - \mu^2)$ and dividing by $\mu \neq 0$ leads to the following characteristic equation:
 
\be \label{eqcarbeta}
2 (\mu + b) \sinh(\lambda) \sin(\lambda) + \left( \alpha + \dfrac{\beta}{\mu} \right) \lambda \left[
\sin(\lambda) \sinh(\lambda \xi) \sinh(\lambda(1 - \xi))- \sinh(\lambda) \sin(\lambda \xi) \sin(\lambda(1 - \xi)) 
\right] = 0.
\ee
\erk

\noindent The asymptotic behaviour of the eigenvalues is often useful for the study of stabilization. The following Proposition gives what we call the asymptotic characteristic equation. 

\bprop (Large eigenvalues of the operator $\caA$) \label{Largeeig}\\
Let $\caA$ be the operator defined by (\ref{defDA}) and (\ref{defA}) in Section \ref{sec2}. Denote by $\mu_n$, $n \in \Z$ the eigenvalues of $\caA$ and by $\lambda$ the complex numbers satisfying $\lambda^4 = - \dfrac{b \mu + \mu^2}{a}$. One of these $\lambda$'s is such that $\exists k \in \Z, \arg(\lambda) \in [-\pi/4 + 2 k \pi;\pi/4 + 2 k \pi)$. The others are $- \lambda$ and $\pm i \lambda$ and they are still denoted by $\lambda$. \\
Assume that $\alpha > 0$, $\alpha b < 6$ and $\xi \notin \Q$. \\
Then $\mu$ is a large eigenvalue of $\caA$ if and only if the corresponding $\lambda$'s satisfy the asymptotic characteristic equation:

\be \label{asymp_eqcar}
\sinh(\lambda) \sin(\lambda) = 0.
\ee

\noindent The algebraic multiplicity of $\mu$ as a large eigenvalue of $\caA$ is exactly $1$ and its geometric multiplicity is $1$.

\eprop  

\proof We start from the characteristic equation (\ref{eqcar1}) in which the trigonometric functions are replaced by:

\be \left \{ \begin{array}{ll}
\sin(\lambda) = \dfrac{e^{\Im(\lambda)}}{2i} \left[ - e^{- i \Re(\lambda)} + e^{- 2 \Im(\lambda)} e^{i \Re(\lambda)} \right] \\
\sinh(\lambda) = \dfrac{e^{\Re(\lambda)}}{2} \left[ e^{i \Im(\lambda)} - e^{- 2 \Re(\lambda)} e^{-i \Im(\lambda)} \right]. 
\end{array} \right.
\ee

\noindent Denoting by  

\be \left \{ \begin{array}{ll}
S(\lambda) := - e^{- i \Re(\lambda)} + e^{- 2 \Im(\lambda)} e^{i \Re(\lambda)} \\
Sh(\lambda) := e^{i \Im(\lambda)} - e^{- 2 \Re(\lambda)} e^{-i \Im(\lambda)} 
\end{array}  \right.
\ee

\noindent the characteristic equation (\ref{eqcar1}) is equivalent to:

\be 
\dfrac{e^{\Im(\lambda)+ \Re(\lambda)}}{4i} \left[2 a \lambda^3 S(\lambda) Sh(\lambda) - \alpha \mu [S(\lambda) Sh(\lambda \xi) Sh(\lambda(\xi - 1)) - Sh(\lambda) S(\lambda \xi) S(\lambda(\xi - 1))] \right] = 0 
\ee

\noindent and, since $\lambda = 0$ is excluded here, it is also:

\be \label{halpha} 
S(\lambda) Sh(\lambda) - \alpha \dfrac{\mu}{2 a \lambda^3} [S(\lambda) Sh(\lambda \xi) Sh(\lambda(\xi - 1)) - Sh(\lambda) S(\lambda \xi) S(\lambda(\xi - 1))]= 0 
\ee 

\noindent The expression $[S(\lambda) Sh(\lambda \xi) Sh(\lambda(\xi - 1)) - Sh(\lambda) S(\lambda \xi) S(\lambda(\xi - 1))]$ is bounded with respect to $\lambda$, if $\Re(\lambda)$ and $\Im(\lambda)$ tend to $+ \infty$ and the quotient $\dfrac{\mu}{2 a \lambda^3}$ tends to zero if $|\lambda|$ tends to $+\infty$. The other cases ($\Re(\lambda)$ and $\Im(\lambda)$ tend to $- \infty$, $\Re(\lambda)$ tends to $\pm \infty$ while $\Im(\lambda)$ tends to $\mp \infty$) are treated similarly. Hence (\ref{asymp_eqcar}). \\
\\
The multiplicity of the large eigenvalues follows from the fact that the roots of $\sin(\lambda) = 0$ are simple as well as those of $\sinh(\lambda)=0$.  The result is proved using Rouché's Theorem. We follow the proof of Lemma 2.4 of \cite{Farah}. \\
We denote by $h_{\alpha} (\lambda)$ the left-hand side of (\ref{halpha}) and define, for $N$ large enough, the curves:

\be \left \{ \begin{array}{ll}
\Gamma_{\pm,n}^{S} := \left \{ z / |z  \pm n \pi| = \dfrac{C_0}{n}\right \}  \\
\Gamma_{\pm,n}^{Sh} := \left \{ z / |z  \pm i n \pi| = \dfrac{C'_0}{n}\right \} 
\end{array}  \right.
\ee   

\noindent Our aim is to choose $C_0$ such that, by Rouché's Theorem, $h_{\alpha}$ has the same roots as $\sinh$ (resp. as $\sin$) inside the curve $\Gamma_{\pm,n}^{Sh}$ (resp. $\Gamma_{\pm,n}^{S}$) for every $n > N$ where $N$ is large enough. \\
The proof is written for $z \in \Gamma_{+,n}^{Sh}$. The rest is analogous. \\
\\
\noindent The first step is to show that, if $z \in \Gamma_{+,n}^{Sh}$, then $|Sh(z)| \geq \dfrac{C_0}{|z|}$. \\
Writing $z = i n \pi + \rho_n e^{i \theta}$ with $\rho_n = \dfrac{C_0}{n}$ and using trigonometric formulae lead to 

\be \begin{array}{ll}
|\sinh(z)|^2 &= |\sinh(i n \pi + \rho_n e^{i \theta})|^2 = \sinh^2(\rho_n \cos(\theta)) \cos^2(\rho_n \cos(\theta)) + \cosh^2(\rho_n \cos(\theta)) \sin^2(\rho_n \cos(\theta)) \\
&= \rho_n^2 + o(1). 
\end{array} 
\ee   

\noindent Now $|z|^2 = \rho_n^2 \cos^2(\theta) +(n \pi + \rho_n \sin(\theta))^2 \leq (n \pi + \rho_n \sin(\theta))^2$ and

\be 
\dfrac{C_0^2}{|z|^2} \leq \dfrac{C_0^2}{(n \pi + \rho_n \sin(\theta))^2} = \dfrac{C_0^2}{n^2 \pi^2} + o \left( \dfrac{1}{n^2} \right) = \dfrac{C_0^2}{n^2} + o(1) = |\sinh(z)|^2.  
\ee

\noindent Since $|Sh(z)| = 2 e^{- \Re(z)} |\sinh(z)|$ and $\Re(z) = \rho_n \cos(\theta)$ tends to zero when $n$ tends to $+\infty$ 

\be  
\exists N_1, n \geq N_1 \Rightarrow 2 e^{-\Re(z)} \geq 1 \; \mbox{and} \; |Sh(z)| \geq |\sinh(z)| \geq \dfrac{C_0}{|z|}. 
\ee

\noindent The second step is to show that, if $z \in \Gamma_{+,n}^{Sh}$, then $|h_{\alpha}(z) - S(z) Sh(z)| \leq |S(z) Sh(z)|$. \\ 
\\
By definition of $h_{\alpha}(z)$, it holds:

\be 
|h_{\alpha}(z) - S(z) Sh(z)| \leq \dfrac{1}{|z|} \cdot \left| \dfrac{\alpha \mu}{2 a \lambda^2} \right| \cdot |S(z)| \cdot \left| Sh(\xi z)Sh((\xi-1) z) - Sh(z) \dfrac{S(\xi z)S((\xi-1) z)}{S(z)} \right|.
\ee

\noindent If $z \in \Gamma_{+,n}^{Sh}$, $\Im(z) = n \pi + \dfrac{C_0}{n} \sin(\theta) \geq n \pi - \dfrac{C_0}{n}$. \\
The expression $\dfrac{\alpha \mu}{2 a z^2} \left( Sh(\xi z)Sh((\xi-1) z) - Sh(z) \dfrac{S(\xi z)S((\xi-1) z)}{S(z)}\right)$ is bounded with respect to $\Im(z)$. Let us denote by $C_0$ the real number such that 

\be 
\exists C_1, \forall z \in \Gamma_{+,n}^{Sh}, \Im(z) \geq C_1 \Rightarrow \left| \dfrac{\alpha \mu}{2 a \lambda^2} \right| \cdot \left| Sh(\xi z)Sh((\xi-1) z) - Sh(z) \dfrac{S(\xi z)S((\xi-1) z)}{S(z)} \right| \leq C_0.
\ee

\noindent At last, there exists $N_2$ such that $N_2 \pi - \dfrac{C_0}{N_2} > C_1$. \\
If $N \geq \max \{N_1,N_2 \}$, $z \in \Gamma_{+,n}^{Sh} \Rightarrow |h_{\alpha}(z) - S(z) Sh(z)| \leq |S(z) Sh(z)|$. \fin   

\brk The large roots of the characteristic equation (\ref{eqcar}) are close to the roots of (\ref{asymp_eqcar}), that is to say, either $(n \pi)$ or $(i n \pi)$ with $n \in \Z^*$. The large eigenvalues $\mu$ of the operator $\caA$ with $\alpha > 0$ are thus near the eigenvalues of the operator $\caA_0$ given in the first part of Theorem \ref{Theig}.
\erk

\brk (Case $\beta>0$) \label{rklarge} The latter Proposition still holds for Problem (\ref{eq1bis})-(\ref{ICbis}) with the additional term $\beta \delta_{\xi}$.
\erk

\noindent Another natural question which arises in this context is that of the continuity of the eigenvalues with respect to the parameter $\alpha$. 

\bprop (Continuity of the eigenvalues and eigenvectors of the operator $\caA$ with respect to $\alpha$) \label{Conteig}\\
Let $\caA$ be the operator defined by (\ref{defDA}) and (\ref{defA}) in Section \ref{sec2}. Denote by $\mu$ an eigenvalue of $\caA$ and by $\lambda$ the complex numbers satisfying $\lambda^4 = - \dfrac{b \mu + \mu^2}{a}$. One of these $\lambda$'s is such that $\exists k \in \Z, \arg(\lambda) \in [-\pi/4 + 2 k \pi;\pi/4 + 2 k \pi)$. The others are $- \lambda$ and $\pm i \lambda$ and they are still denoted by $\lambda$. \\
Assume that $\alpha > 0$, $\alpha b < 6$ and $\xi \notin \Q$. \\
Then $\mu$ depends continuously on the parameter $\alpha$ and the eigenvector $\Phi(\mu, \cdot)$ converges pointwise to the eigenvector $\Phi_n^{\pm}$ (up to a multiplicative constant) when $\alpha$ tends to zero.
\eprop

\proof To state the continuity of $\mu$ with respect to $\alpha$, we follow the proof of Remark $1$ of \cite{Farah}. We denote by $g_{\alpha} (\lambda)$ the left-hand side of (\ref{eqcar}).\\ 
For a fixed $\alpha$, denote by $\lambda_0$ a root of $\caA$. Since $\lambda_0$ is isolated, there exists $\rho > 0$ such that

\be 
g_{\alpha}(z) \neq 0, \forall z \in \C, \; \mbox{such that} \; 0 < |z - \lambda_0 | < \rho.
\ee

\noindent Now, $g_{\alpha}$ is a continuous function. Indeed the functions $\lambda \mapsto \lambda^3$, $\lambda \mapsto \sin(\lambda)$ and $\lambda \mapsto \sinh(\lambda)$ are continuous functions on $\C$. As for $\mu$, it is the root of the second degree equation:

\be
\mu^2 + b \mu + a \lambda^4 = 0.
\ee 

\noindent Thus it is a continuous function of the coefficients of this equation i.e. a continuous function of the variable $\lambda$. \\
Setting $D = \{ z \in \C, \; \mbox{such that} \; |z - \lambda_0 | = \rho \}$, the continuity of $g_{\alpha}$ implies that there exists a positive real number $\kappa$ such that $|g_{\alpha}(z)| \geq \kappa, \forall z \in D$. \\
For a fixed positive real number $\epsilon_0$, we consider the mapping of two variables

\be
H : [0, \epsilon_0] \times D \rightarrow \C : (\epsilon, z) \mapsto g_{\alpha + \epsilon}(z) - g_{\alpha}(z). 
\ee
\noindent Since it is a uniformly continuous function and since $H(0; z) = 0$ for all $z \in D$, we deduce the existence of a positive real number $\delta$ such that
\be
|H(\epsilon, z)| < \kappa, \forall (\epsilon, z) \in [0;\delta] \times D. 
\ee

\noindent The last two estimates imply that
\be
|g_{\alpha + \epsilon}(z) - g_{\alpha}(z)| < |g_{\alpha}(z)|, \forall (\epsilon, z) \in [0;\delta] \times D. 
\ee

\noindent Hence Rouché's theorem allows to conclude that $g_{\alpha + \epsilon}$ has the same number of roots as $g_{\alpha}$ for all $\epsilon \in [0;\delta]$. Thus the root $\lambda$ of $g_{\alpha}$ is a continuous function of $\alpha$ and, by composition, $\mu$ is also a continuous function of $\alpha$.   \\
\\
\noindent To finish with the proof, note that, when $\alpha$ tends to zero, $\alpha b < 6$ then $\mu + b \neq 0$ and the characteristic equation (\ref{eqcar}) becomes: $\sin(\lambda) \sinh(\lambda)=0$. Thus $\lambda$ tends either to $(n \pi)$ or to $(i n \pi)$ with $n \in \Z^*$. \\
If it tends to $(n \pi)$ for example, $\sin(\lambda)$ tends to zero and the eigenfunction tends to:

\be \label{limeigfunc1}
\lim_{\lambda \rightarrow n \pi} \phi(\lambda, x):= \dfrac{1}{n^2 \pi^2} \cdot e^{n \pi (\xi - 2)} \sinh(n \pi) \sin(n \pi (1-\xi)) \sin(n \pi x). 
\ee

\noindent Indeed the other terms tend to zero since $\sin(\lambda)$ tends to zero and the following two expressions are bounded with $\lambda$:

\be \left \{ \begin{array}{ll}
\dfrac{1}{|\lambda|^2} \cdot e^{(\xi - 2) |\lambda|} \sinh(\lambda) [\sin(\lambda(x - \xi)) - \sinh(\lambda(x - \xi))] \cdot H(x, \xi) \\
\dfrac{1}{|\lambda|^2} \cdot e^{(\xi - 2) |\lambda|} \sinh(\lambda(1-\xi)) \sinh(\lambda x). 
\end{array} \right.
\ee

\noindent And if $\lambda$ tends to $(i n \pi)$, $\sinh(\lambda)$ tends to zero and the eigenfunction tends to:

\be \label{limeigfunc2}
\begin{split}
\lim_{\lambda \rightarrow i n \pi} \phi(\lambda, x):= & \dfrac{1}{n^2 \pi^2} \cdot e^{n \pi (\xi - 2)} \sin(i n \pi) \sinh(i n \pi (1-\xi)) \sinh(i n \pi x) \\
=& - i \dfrac{1}{n^2 \pi^2} \cdot e^{n \pi (\xi - 2)} \sinh(n \pi) \sin(n \pi (1 - \xi) \sin(n \pi x).
\end{split}  
\ee

\noindent Indeed the other terms tend to zero since $\sinh(\lambda)$ tends to zero and the following two expressions are bounded with $\lambda$:

\be \left \{ \begin{array}{ll}
\dfrac{1}{|\lambda|^2} \cdot e^{(\xi - 2) |\lambda|} \sin(\lambda) [\sin(\lambda(x - \xi)) - \sinh(\lambda(x - \xi))] \cdot H(x, \xi) \\
\dfrac{1}{|\lambda|^2} \cdot e^{(\xi - 2) |\lambda|} \sin(\lambda(1-\xi)) \sin(\lambda x). 
\end{array}  \right.
\ee

\noindent Hence the announced result concerning the eigenfunctions. \fin


\section{Localization of the eigenvalues of the operator $\caA$ for small values of $\alpha$} \label{sec5}


The aim of this Section is to determine the localization of the eigenvalues of the operator $\caA$ for small values of $\alpha$ compared with the case $\alpha =0$ already studied in \cite{AmDiZer} and recalled in Theorem \ref{Theig}. 
\bthe (Localization of the eigenvalues of the operator $\caA$ for small values of $\alpha$) \label{Thloceig}\\
Let $\caA$ be the operator defined by (\ref{defDA}) and (\ref{defA}) in Section \ref{sec2} with $\xi \notin \Q$. Denote by $\mu$ an eigenvalue of $\caA$. \\
\\
\noindent The case $\alpha=0$ is already known (cf. Theorem 
\ref{Theig}). If $b \leq 2 \sqrt{a} \pi^2$, then 

\be
\sup_{\mu \in \sigma(\caA)} \Re(\mu) = - \dfrac{b}{2}. 
\ee

\noindent If $b > 2 \sqrt{a} \pi^2$, then 

\be \label{maxRe}
\sup_{\mu \in \sigma(\caA)} \Re(\mu) = \dfrac{1}{2} (-b + \sqrt{b^2 - 4 a \pi^4}).
\ee 

\noindent Now, for any $\epsilon > 0$, there exists $\alpha_0$ such that, if $0< \alpha < \alpha_0$, it holds: if $b \leq 2 \sqrt{a} \pi^2$, then 

\be
\Re(\mu) \leq - \dfrac{b}{2} + \epsilon. 
\ee

\noindent If $b > 2 \sqrt{a} \pi^2$, then 

\be \label{maxRe}
\Re(\mu) \leq \dfrac{1}{2} (-b + \sqrt{b^2 - 4 a \pi^4}) + \epsilon.
\ee 

\ethe 

\proof Since the study is restricted here to the case of the small non-vanishing values of $\alpha$, the eigenvalues $\mu$ of the operator $\caA$ are known to be such that the corresponding $\lambda$'s satisfy $\lambda \notin \{ k \pi, k \in \Z \}$, $\lambda \notin \{ i k \pi, k \in \Z \}$ and the characteristic equation:

\be \label{eqcar_recall}
2 a \lambda^3 \sinh(\lambda) \sin(\lambda) - \alpha \mu \left[
\sin(\lambda) \sinh(\lambda \xi) \sinh(\lambda(1 - \xi))- \sinh(\lambda) \sin(\lambda \xi) \sin(\lambda(1 - \xi)) 
\right] = 0
\ee 

\noindent (cf. Proof of Theorem \ref{Theig}). \\
\\
\noindent The first remark is that this equation is invariant under the transform $\lambda \mapsto i \lambda$ so it is enough to consider $\lambda \in \{ z \in \C - \{ k \pi, k \in \Z \},  \exists l \in \Z, \arg(z) \in [-\pi/4 + 2 l \pi;\pi/4 + 2 l \pi)$. \\
\\
In order to compute the first degree Mac Laurin polynomial for the function $\lambda(\alpha)$ i.e. to find $\lambda_1$ such that $\lambda = n \pi + \lambda_1 \alpha + o(\alpha)$ around the value $\alpha = 0$, we write all the first degree Mac Laurin polynomials for the functions involved in equation (\ref{eqcar_recall}) as functions of $\lambda$ around $\lambda = n \pi$ ($n \in \N^*$):

\be
\left \{ \begin{array}{lll}
\lambda^3 = n^3 \pi^3 + (3 n^2 \pi^2 \lambda_1) \alpha + o(\alpha) \\
\sinh(\lambda) = \sinh(n \pi) + \cosh(n \pi) \lambda_1 \alpha + o(\alpha) \\
\sin(\lambda) = (-1)^n \lambda_1 \alpha + o(\alpha) 
\end{array} \right.
\ee

\noindent Hence $2 a \lambda^3 \sinh(\lambda) \sin(\lambda) = [2 a (-1)^n n^3 \pi^3 \sinh(n \pi) \lambda_1] \alpha + o(\alpha)$. Now, $\exists n_0$ such that $2 a n_0^2 \pi^2 < b < 2 a (n_0 + 1)^2 \pi^2$ and for $n > n_0$ : 

\be
\mu_n^{+,-} = \dfrac{1}{2} \left[ - b \pm i \sqrt{4 a n^4 \pi^4 - b^2} \right] + o(1). 
\ee

\noindent Inserting all these results into the characteristic equation (\ref{eqcar_recall}) leads, for $n$ such that $n > n_0$, to:

\be \label{lambda1}
\lambda_1 = \dfrac{[- b \pm i \sqrt{4 a n^4 \pi^4 - b^2}] \sin(n \pi \xi) \sin(n \pi (\xi - 1))}{4 (-1)^n a n^3 \pi^3} = \dfrac{[- b \pm i \sqrt{4 a n^4 \pi^4 - b^2}] (-1)^n \sin^2(2 n \pi \xi)}{8 (-1)^n a n^3 \pi^3}.
\ee 

\noindent At last, $\lambda = n \pi + \lambda_1 \alpha + o(\alpha)$ is put into the expression of $\mu_n^{\pm}$, which leads after some calculations to: 

\be \label{mun}
\left \{ \begin{array}{ll}
\mu_n^{\pm} = \dfrac{1}{2} [- b \pm i \sqrt{4 a n^4 \pi^4 - b^2}]  \pm i \dfrac{\sin^2(2 n \pi \xi)}{2 \sqrt{4 a n^4 \pi^4 - b^2}} [- b \pm i \sqrt{4 a n^4 \pi^4 - b^2}] \cdot \alpha + o(\alpha), \forall n > n_0 \\
\mu_n^{\pm} = \dfrac{1}{2} [- b \pm \sqrt{b^2 - 4 a n^4 \pi^4}]  \mp  \dfrac{\sin^2(2 n \pi \xi)}{2 \sqrt{b^2 - 4 a n^4 \pi^4}} [- b \pm \sqrt{b^2 - 4 a n^4 \pi^4}] \cdot \alpha + o(\alpha), \forall n, 1 \leq n \leq n_0.
\end{array} \right.
\ee

\noindent The real part of $\mu_n^{\pm}$ is:

\be \label{Remun}
\left \{ \begin{array}{ll}
\Re(\mu_n^{\pm}) = - \dfrac{b}{2} - \dfrac{1}{2} \sin^2(2 n \pi \xi) \cdot \alpha + o(\alpha), \forall n > n_0 \\
\Re(\mu_n^{\pm}) = \dfrac{1}{2} [- b \pm \sqrt{b^2 - 4 a n^4 \pi^4}]  \pm  \dfrac{\sin^2(2 n \pi \xi)}{2 \sqrt{b^2 - 4 a n^4 \pi^4}} [b \mp \sqrt{b^2 - 4 a n^4 \pi^4}] \cdot \alpha + o(\alpha), \forall n, 1 \leq n \leq n_0.
\end{array} \right.
\ee 

\noindent The eigenvalues of $\caA$ for small positive values of $\alpha$ and large values of $n$ are on the left of those we have for $\alpha =0$ (whose real part is $-b/2$). \\
\\
As for the case of small values of $n$ ($1 \leq n < n_0$), since $[b \mp \sqrt{b^2 - 4 a n^4 \pi^4}] > 0$ and $\sqrt{b^2 - 4 a n^4 \pi^4} > 0$, the eigenvalues $\mu_n^+$ of $\caA$ for small positive values of $\alpha$ are on the right of those we have for $\alpha =0$ and the eigenvalues $\mu_n^-$ of $\caA$ for small positive values of $\alpha$ are on the left of those we have for $\alpha =0$. \\
\\
\noindent Moreover, $0 \leq \dfrac{1}{2} \sin^2(2 n \pi \xi) \cdot \alpha \leq \dfrac{\alpha}{2}$ and $\left| \dfrac{\sin^2(2 n \pi \xi)}{2 \sqrt{b^2 - 4 a n^4 \pi^4}} [b \mp \sqrt{b^2 - 4 a n^4 \pi^4}] \right|$ is bounded for $n \leq n_0$. \\  
\\
Now, it holds, for all $1 \leq n \leq n_0$:

\be \label{ordsmall}
- b < \dfrac{1}{2} (-b - \sqrt{b^2 - 4 a \pi^4}) < \dfrac{1}{2} (-b - \sqrt{b^2 - 4 a n^4 \pi^4}) < - \dfrac{b}{2} < \dfrac{1}{2} (-b + \sqrt{b^2 - 4 a n^4 \pi^4}) < \dfrac{1}{2} (-b + \sqrt{b^2 - 4 a \pi^4}) < 0.
\ee

\noindent Thus, for all $n \geq 1$ and $\alpha$ sufficiently small:

\be \label{maxRe}
\Re(\mu_n^{\pm}) < \dfrac{1}{2} (-b + \sqrt{b^2 - 4 a \pi^4}).
\ee   

\noindent Note that, if $\xi$ belonged to $\Q$, there would exist $(m,n) \in (\N^*)^2$ such that $\xi = \dfrac{m}{n}$ and $\sin(2 n \pi \xi)$ would vanish (as well as $\sin(2 k n \pi \xi)$ for $k \in \Z$). More calculations would be required to get the third (at least) degree Mac Laurin polynomials of the functions involved in the characteristic equation since the second degree term for $\lambda$ and for $\mu_n$ also vanishes. \fin 

\brk (Case $\beta > 0$) \label{rktaylor} The first-degree Taylor polynomial linear approximation of $\Re(\mu_n^{\pm})$ if we add the term $\beta \delta_{\xi}$ to get Problem (\ref{eq1bis})-(\ref{ICbis}) contains the additional term:

\be 
\dfrac{\pm 2 \sin^2(2 n \pi \xi)}{\sqrt{b^2 - 4 a n^4 \pi^4}} \cdot \beta.
\ee

\noindent Note that the real part of the eigenvalues $\mu_n^{+}$ is larger than that of the same eigenvalues with $\beta = 0$. This confirms what we said in the introduction: the stiffness term does not improve the decreasing of the energy. If the SMA cables restrain the vibrations of girder bridges, it only comes from the damper term $\alpha u_t \delta_{\xi}$. \\
All that is true if $\sin^2(2 n \pi \xi)$ does not vanish, which happens for some values of $\xi$. That is why this result does not contradict the experiments of \cite{LiuFu}.  
\erk


\section{Energy decreasing} \label{sec7}

\begin{definition} (Functions $\Psi_n^{\pm}$) \\
Consider $\mu$ such that $\lambda$ defined by $\lambda^4 = - \dfrac{b \mu + \mu^2}{a}$ satisfies the characteristic equation (\ref{eqcar}).  It is already known that $\mu$ is an eigenvalue of the operator $\caA$ and that one of these $\lambda$'s is such that $\exists k \in \Z, \arg(\lambda) \in [-\pi/4 + 2 k \pi;\pi/4 + 2 k \pi)$. The others are $- \lambda$ and $\pm i \lambda$. Consider this $\lambda$. If its modulus $|\lambda|$ tends to $+ \infty$, then there exists $n \in \Z^*$ such that $\lambda$ tends to $n \pi$ (respectively to $i n \pi$). \\
This integer $n$ depends on $\mu$. \\
The function $\Psi_n^{\pm}$ is defined on $(0;1)$ by

\be \label{Psi}
\Psi_n^{\pm}:= \dfrac{1}{n^2 \pi^2} e^{n \pi(\xi - 2)} \sinh(n \pi) \sin(n \pi(1 - \xi)) \sin(n \pi x)(1, \mu^{\pm})
\ee 

\noindent (respectively by
\be \label{Psi2}
\Psi_n^{\pm}:= -i \dfrac{1}{n^2 \pi^2} e^{(\xi - 2) |\lambda|} \sin(n \pi) \sin(n \pi x)(1, \mu^{\pm})).
\ee 

\end{definition}


\begin{theorem} \label{rieszbasis} (Riesz basis for the operator $\mathcal{A}$) \\ 
Let $\Phi(\mu, \cdot)$ still be defined as in Theorem \ref{Theig} and $\Psi_n^{\pm}$ given by the above definition. There exists $n_0 \geq 1$, such that

\be 
\sum_{|n| \geq n_0} \| \Phi(\mu, \cdot) - \Psi_n^{\pm} \|_{\caH}^2 < \infty.
\ee

\noindent Thus, the root eigenvectors of $\caA$ form a Riesz basis of $\caH$. 
\end{theorem}

\proof First the $\Psi_n^{\pm}$ form an orthogonal basis (see \cite{AmDiZer}). \\
\\
\noindent The inner product in $\caH$ has two terms. Since the eigenvector $\Phi(\mu, \cdot)$ is defined as $\Phi(\mu, \cdot):=(\phi(\lambda, \cdot)(1,\mu)$ and since $\mu$ has the same behaviour as $\lambda^2$ for large values of $\lambda$, it is enough to consider the first term in the inner product. \\
Now the eigenfunction defined by (\ref{eigfunc}) is made of three terms. Let us start with the third one. The second one is analogous with easier calculations. That is why we do not give details for this second term. 

\noindent Let us prove that, if $\lambda \rightarrow n \pi$ ($n \in \N^*$), then

\be \label{term3}
\begin{split}
\int_0^1 \left| e^{(\xi - 2) |\lambda|} \sin(\lambda) \sinh(\lambda(1-\xi)) \sinh(\lambda x) + e^{(\xi - 2) n \pi} \sin(n \pi) \sinh(n \pi(1-\xi)) \sinh(n \pi x) \right|^2 dx = O \left( \dfrac{1}{n^2} \right).  
\end{split}
\ee

\noindent The integrand is of the form $|A \cdot B - A_0 \cdot B_0|^2 \leq 2 |A - A_0|^2 \cdot |B|^2 + 2 |A_0|^2 \cdot |B - B_0|^2$ with

\be 
|A_0|^2 : = |e^{(\xi - 2) n \pi} \sin(n \pi) \sinh(n \pi(1-\xi))|^2 = O \left(e^{-2 n \pi} \right). 
\ee

\noindent Since, for $x \in (0;1)$, $\sinh(\lambda x) = O \left( e^{2 |\Re(\lambda)|}\right) = O \left( e^{2 n \pi}\right)$ (cf. proof of Proposition \ref{Largeeig} about the large eigenvalues), it holds:

\be 
\int_0^1 |B|^2 : = \int_0^1 \left| \sinh(\lambda x) \right|^2 dx = O \left( e^{2 n \pi} \right). 
\ee

\noindent Thus it is enough to show the following two estimates:

\be \label{AminusA0}
\begin{split}
|A - A_0|^2  := \left| e^{(\xi - 2) |\lambda|} \sin(\lambda) \sinh(\lambda(1-\xi)) - e^{(\xi - 2) n \pi} \sin(n \pi) \sinh(n \pi(1-\xi)) \right|^2 = O \left( \dfrac{1}{n^2} e^{-2 n \pi} \right)  
\end{split}
\ee

\be \label{BminusB0}
\int_0^1 |B - B_0|^2 := \int_0^1 \left| \sinh(\lambda x) - \sinh(n \pi x) \right|^2 dx = O \left( \dfrac{1}{n^2} e^{2 n \pi} \right). 
\ee

\noindent Coming back to the definition of the trigonometric functions, it holds:
\be 
\sin(\lambda) - \sin(n \pi) = \dfrac{1}{2i} \left[ (e^{i \Re(\lambda)} e^{- \Im(\lambda)} - e^{i n \pi}) - (e^{-i \Re(\lambda)} e^{\Im(\lambda)} - e^{-i n \pi}) \right]. 
\ee

\noindent Then, using the mean-value Theorem as well as the asymptotic behaviour of the difference $(\lambda - n \pi)$ obtained in the proof of Proposition \ref{Largeeig} about the large eigenvalues: 

\be \begin{array}{llll}
\left| e^{i \Re(\lambda)} e^{- \Im(\lambda)} - e^{i n \pi} \right|^2 = \left| \left( e^{i \Re(\lambda)} - e^{i n \pi} \right) + \left( e^{- \Im(\lambda)} - 1 \right) \cdot e^{i \Re(\lambda)} \right|^2 \\
\leq 2  \left( e^{- \Im(\lambda)} - 1 \right)^2 + 2 \left| \cos(\Re(\lambda)x) - \cos(n \pi) + i \left( \sin(\Re(\lambda)) - \sin(n \pi) \right) \right|^2 \\
\leq 2  \left( e^{- \Im(\lambda)} - 1 \right)^2 + 2 \left( \cos(\Re(\lambda)) - \cos(n \pi) \right)^2 + 2 \left( \sin(\Re(\lambda)) - \sin(n \pi) \right)^2 \\
\leq  \dfrac{C_0^2}{n^2} \exp \left(2 \left( \dfrac{C_0}{n} + o(1/n) \right) \right) + 4 \left( \Re(\lambda) - n \pi \right)^2 
= O \left( \dfrac{1}{n^2} \right). 
\end{array}
\ee

\noindent Now (\ref{AminusA0}) is of the form $|E \cdot F - E_0 \cdot F_0|^2 \leq 2 |E - E_0|^2 \cdot |F|^2 + 2 |E_0|^2 \cdot |F - F_0|^2$ with 

\be \left \{ \begin{array}{llll}
|F - F_0|^2:= \left| e^{i \Re(\lambda)} e^{- \Im(\lambda)} - e^{i n \pi} \right|^2 = O \left( \dfrac{1}{n^2} \right) \\
|E_0|^2 := \left| e^{(\xi - 2) n \pi} \sinh(n \pi(1-\xi)) \right|^2 = O \left( e^{-2 n \pi} \right) \\
|E - E_0|^2:= \left| e^{(\xi - 2) |\lambda|} \sinh(\lambda(1-\xi)) - e^{(\xi - 2) n \pi} \sinh(n \pi(1-\xi)) \right|^2 = O \left( e^{-2 n \pi} \right) \\
|\sin(\lambda)|^2 = \left| \sin \left( n \pi + O \left( \dfrac{1}{n} \right) \right) \right|^2 = \left| (-1)^n \sin \left(O \left( \dfrac{1}{n} \right) \right) \right|^2 = O \left( \dfrac{1}{n^2} \right).
\end{array} \right.
\ee

\noindent Hence (\ref{AminusA0}). Analogously, it holds:

\be 
\sinh(\lambda x) - \sinh(n \pi x) = \dfrac{1}{2} \left[ (e^{\Re(\lambda) x} e^{i \Im(\lambda) x} - e^{n \pi x}) - (e^{-\Re(\lambda) x} e^{- i \Im(\lambda) x} - e^{-n \pi x}) \right]. 
\ee

\noindent Using once more the mean-value Theorem as well as the asymptotic behaviour of the difference $(\lambda - n \pi)$ obtained in the proof of Proposition \ref{Largeeig} about the large eigenvalues:  

\be \begin{array}{llll}
\left| e^{\Re(\lambda) x} e^{i \Im(\lambda) x} - e^{n \pi x} \right|^2 = \left| \left( e^{\Re(\lambda) x} - e^{n \pi x} \right) + \left( e^{i \Im(\lambda) x} - 1 \right) \cdot e^{\Re(\lambda) x} \right|^2 \\
\leq 2 \left( e^{\Re(\lambda) x} - e^{n \pi x} \right)^2 + 2 e^{2 \Re(\lambda) x} \left[ \left(\cos(\Im(\lambda)x) - 1 \right)^2 + \sin^2(\Im(\lambda)x) \right]  \\
\leq 4 \dfrac{C_0^2}{n^2} \exp \left(2 \left( n \pi + \dfrac{C_0}{n} + o(1/n) \right) \right) \\
= O \left( \dfrac{1}{n^2} e^{2 n \pi} \right). 
\end{array}
\ee

\noindent Hence (\ref{BminusB0}). \\
\\
The first term of the eigenfunction has to be considered. \\
Since $|\sin(\lambda)|^2 = O \left( \dfrac{1}{n^2} \right)$ and $\left| e^{(\xi - 2) |\lambda|} \sinh(\lambda) \sinh(\lambda (x- \xi)) H(x,\xi) \right|^2 = O(1)$ for $x > \xi$, it holds, if $\lambda \rightarrow n \pi$ ($n \in \N^*$):

\be \label{term1}
\begin{split}
\int_0^1 \left| e^{(\xi - 2) |\lambda|} \sin(\lambda) \sinh(\lambda) \sinh(\lambda (x- \xi)) H(x,\xi) \right|^2 dx = O \left( \dfrac{1}{n^2} \right).  
\end{split}
\ee
 
\noindent This ends the proof. \fin \\
\\
\noindent Finally, the energy is proved to decay exponentially and the localization of the eigenvalues for small values of $\alpha$ leads to a lower bound for the decay rate in that case. This bound is the optimal decay rate obtained for $\alpha = 0$ in \cite{AmDiZer} which means that the decay rate for small values of $\alpha > 0$ is equal to that with $\alpha =0$ (i.e. without SMA cables). SMA cables do not restrain the vibrations of girder bridges if $\alpha$ is close to zero. This is not in contradiction with the experiments of \cite{LiuFu} but it means that the simple fact of adding SMA cables may not be enough. 

\bthe \label{expo}(Exponential stability and decay rate for small values of $\alpha$)   \\
System (\ref{eq1})-(\ref{IC}) (presented in the introduction) is exponentially stable and, for any $\epsilon > 0$, there exists $\alpha_0$ such that for any $0 < \alpha < \alpha_0$, for any solution $U$ of (\ref{abstr}) with initial data $U_0 \in D(\caA)$, there exist constants $C>0$ and $\omega_0(\epsilon) <0$ depending on $a$, $b$ and $\alpha$ such that:
$$E(t) \leq C e^{2 \omega_0(\epsilon) t} \| U_0 \|^2_{D(\caA)}, \forall t > 0,$$ where $E(t)$ is defined by (\ref{defE}) and 

\be
\left \{ \begin{array}{ll}
\omega_0(\epsilon) : = -\dfrac{b}{2} + \epsilon, \; \mbox{if} \; b \leq 2 \sqrt{a} \pi^2  \\ 
\omega_0(\epsilon) : = \dfrac{1}{2} (-b + \sqrt{b^2 - 4 a \pi^4}) + \epsilon, \; \mbox{if} \; b > 2 \sqrt{a} \pi^2.
\end{array} \right.
\ee
 
\ethe

\proof According to Theorem \ref{rieszbasis}, the system of eigenvectors of $\caA$ constitutes a Riesz basis. Consequently, by a standard argument (see the proof of Theorem $2.5$ of \cite{AmDiZer} for example), the optimal energy decay rate is identified with the supremum of the real part of the eigenvalues of $\caA$. \\
Thus the result follows from Theorem \ref{Thloceig}. \fin 

\brk (Case $\alpha = 0$) \\
In \cite{AmDiZer} the systems of eigenvectors of $\caA$ and $\caA_0$ (obtained with $\alpha = 0$) are proved to be quadratically close in $V×L_2(0,1)$ using the explicit expression of the eigenfunctions of both operators. Thus, it follows from Theorem 3 in Appendix D of \cite{PTru}, that the system of eigenvectors of $\caA_0$ constitutes a Riesz basis. Consequently, by a standard argument (see the proof of Theorem $2.5$ in \cite{AmDiZer}), the optimal energy decay rate is identified with the supremum of the real part of the eigenvalues of $\caA_0$. With the notation of Theorem \ref{expo}, it means that $\omega = \omega_0(0)$ when $\alpha = 0$. 
\erk


\section{Conclusion} \label{sec9}

One could have thought that adding the damping term $\alpha \delta_{\xi} u_t$, $\alpha > 0$ would logically increase the decay rate of the energy compared with the case $\alpha = 0$ already studied by Ammari, Dimassi and Zerzeri in \cite{AmDiZer}. \\
\\
In fact this is not obvious. If $\alpha$ is small, we have proved that the decay rate does not change. Is it clear that it will increase with $\alpha$? This is not clear either. \\
\\
As announced in the introduction, let us prove that, if $\alpha=0$, increasing the value of $b$ may not increase the decay rate of the energy. Indeed, in that situation, the energy decreases exponentially and the decay rate is given by the maximum of $\Re{\mu_n^{+}}$ which is   

\be
\Re{\mu_n^{+}} = \dfrac{1}{2} \left( -b + \sqrt{b^2 - 4 a n^4 \pi^4} \right) = \dfrac{-2 a n^4 \pi^4}{b + \sqrt{b^2 - 4 a n^4 \pi^4}}, n = 1,2, \ldots, n_0 \\
\ee

\noindent where $b \in \R^{+*}- \{ 2 \sqrt{a} n^2 \pi^2, n \in \Z^* \}$ and $n_0$ is such that $2 \sqrt{a} n_0^2 \pi^2 < b < 2 \sqrt{a} (n_0 + 1)^2 \pi^2$. \\
Thus

\be
\Re{\mu_n^{+}} \leq \dfrac{-2 a \pi^4}{b + \sqrt{b^2 - 4 a n_0^4 \pi^4}}, n = 1,2, \ldots, n_0 \\
\ee  

\noindent and the last bound tends to zero, if $b$ tends to $+ \infty$. \\
\\
\noindent In conclusion, the decay rate of the energy may be better with $\alpha > 0$ (with $\alpha$ large enough) than with $\alpha = 0$ i.e. Shape Memory Alloy cables may restrain the vibrations of girder bridges effectively. \\
 \\
The experiments of \cite{LiuFu} would be confirmed mathematically if we were able to prove that there exist $\alpha > 0$, $\beta > 0$, $\xi \in (0;1)$ and $\omega_1 (\alpha, \beta, \xi) < \omega_0(0)$, such that all the roots of the characteristic equation (\ref{eqcarbeta}) are such that the corresponding $\mu$'s satisfy $\Re(\mu) \leq \omega_1(\alpha, \beta, \xi)$. We conjecture that this requires a value for $\xi$ which is excluded in our theorems.

\vspace{1cm}

\noindent \bf Acknowledgements. \rm Thanks to S. Nicaise who suggested this problem to me as well as the reflection around the fact that increasing the value of $b$ when $\alpha = 0$ may not increase the decay rate. \\
I would like to dedicate this paper to D. Mercier with whom I have enjoyed doing mathematics for so many years. I wish him a happy retirement!


{}
 
\edc
\begin{thebibliography}{99}

\bibitem{Farah} F. Abdallah, D. Mercier, S. Nicaise, Spectral analysis and exponential or polynomial stability of some indefinite sign damped problems, Evol. Equ. Control Theory \bf{2, No 1} \rm(2013) 1-33.

\bibitem{AkilChiGha} M. Akil, Y. Chitour, M. Ghader, A. Wehbe, Stability and Exact Controllability of a Timoshenko System with Only One Fractional Damping on the Boundary, Asymptotic Analysis \bf{119} \rm(2020) 221-280.


\bibitem{AkilAssaWehbe} M. Akil, I. Issa, A. Wehbe, Energy decay of some boundary coupled systems involving wave Euler-Bernoulli beam with one locally singular fractional Kelvin-Voigt damping, arXiv:2102.12732.



\bibitem{AmDiZer} K. Ammari, M. Dimassi, M. Zerzeri, The rate at which energy decays in a viscously damped hinged Euler-Bernoulli beam, J. Diff. Equ. \bf{257} \rm(2014) 3501-3520.

\bibitem{AmTuc00} K. Ammari, M. Tucsnak, Stabilization of Bernouilli-Euler beams by means of a pointwise feedback force, SIAM J. Control Optim. \bf{39} \rm (No. 4) (2000) 1160-1181.

\bibitem{AmTuc01} K. Ammari, M. Tucsnak, Stabilization of  second order evolution equations by a class of unbounded feedbacks, ESAIM: COCV \bf{6} \rm (2001) 361-386.

\bibitem{AB} W. Arendt, C.J.K. Batty, Tauberian theorems and stability of one-parameter of semi-groups, Trans. Amer. Math. Soc. \bf{305 (2)} \rm (1988) 837-852.

\bibitem{Bor} A.~Borichev and Y.~Tomilov, Optimal polynomial decay of functions and operator semigroups, Math. Ann. \textbf{347(2)} \rm(2010) 455-478.

\bibitem{caszu} C. Castro, E. Zuazua, Exact boundary controllability of two Euler-Bernoulli beams connected by a point mass, Mathematical and Computer Modelling, \bf32\rm (2000) 955-969.

\bibitem{coxzu} S. Cox and E. Zuazua, The rate at which energy decays in a damped string, Partial Differential Equations, \bf19\rm, (1994) 213-243.

\bibitem{deko1} B. Dekoninck, S. Nicaise, Control of networks of Euler-Bernoulli beams, ESAIM : Control, Optimisation and Calculus of Variations, \bf4\rm (1999) 57-81.

\bibitem{Maj} L. Majkut, Eigenvalue based inverse model of beam for structural modification and diagnostics. Part I: Theoretical formulation, Latin American Journal of Solids and Structures (2010) 423-436.

\bibitem{LiuFu} A-R. Liu, C-H. Liu, J-Y. Fu, Y-L. Pi, Y-H. Huang, J-P. Zhang, A Method of Reinforcement and Vibration Reduction of Girder Bridges Using Shape Memory Alloy Cables, Int. J. Struct. Stab. Dyn. \bf{17} \rm (No. 7) (2017) 1750076.

\bibitem{merreg1} D. Mercier, V. Régnier, Spectrum of a network of Euler-Bernoulli beams. J. Math. Anal. Appl., \bf 337/1 \rm (2007) 174-196.

\bibitem{merreg2} D. Mercier, V. Régnier, Control of a network of Euler-Bernoulli beams. J. Math. Anal. Appl., \bf 342 \rm (2008) 874-894.

\bibitem{merreg3} D. Mercier, V. Régnier, Boundary controllability of a chain of serially connected Euler-Bernoulli beams with interior masses. Collect. Math. \bf 60/3 \rm (2009) 307-334.

\bibitem{EECT18} D. Mercier, V. Régnier, Decay rate of the Timoshenko system with one boundary damping, Evol. Equ. Control Theory \bf{8, No 2} \rm(2019) 423-445.

\bibitem{P} A. Pazy, Semigroups of Linear Operators and Applications to Partial Differential Equations, Applied Mathematical Sciences \bf{44} \rm(1983) Springer-Verlag.

\bibitem{PTru} J. P\"oschel, E. Trubowitz, Inverse Spectral Theory, Pure and Applied Mathematics, \bf130\rm, Academic Press, Boston, MA (1987).

\bibitem{young} R. M. Young, An Introduction to Nonharmonic Fourier Series, Academic Press, New York (1980).

\end{thebibliography}
